\documentclass[leqno,12pt]{article}

\usepackage{amssymb,hyperref,amsthm}
\usepackage{euscript}
\usepackage[dvips]{graphicx}
\usepackage{flafter}
\usepackage{pstricks}
\usepackage[latin1]{inputenc}
\usepackage{amsmath}
\usepackage{amsfonts}
\usepackage{pstricks-add}
\usepackage{dsfont}
\usepackage{bm}

\usepackage[refpage]{nomencl}
\makenomenclature

\usepackage{makeidx}
\makeindex

\usepackage{mathrsfs} 

\evensidemargin\oddsidemargin

\begin{document}

\newcommand{\eqnsection}{
\renewcommand{\theequation}{\thesection.\arabic{equation}}
   \makeatletter
   \csname  @addtoreset\endcsname{equation}{section}
   \makeatother}
\eqnsection

\def\r{{\mathbb R}}
\def\e{{\mathbb E}}
\def\p{{\mathbb P}}
\def\P{{\bf P}}
\def\E{{\bf E}}
\def\Q{{\bf Q}}
\def\z{{\mathbb Z}}
\def\N{{\mathbb N}}
\def\T{{\mathbb T}}
\def\G{{\mathbb G}}
\def\L{{\mathbb L}}
\def\1{{\mathds{1}}}
\def\deg{\chi}

\def\d{\mathtt{d}}
\def\ttheta{{\bm \theta}}
\def\t{{\bf{t}}}
\def\a{{\bf{a}}}
\def\deg{\chi}
\def\B{\mathfrak{B}}

\def\M{{\mathbb{M}}}
\def\ee{\mathrm{e}}
\def\d{\, \mathrm{d}}
\def\S{\mathscr{S}}
\def\bs{{\tt bs}}
\def\bbeta{{\bm \beta}}
\def\ttheta{{\bm \theta}}

\newtheorem{theorem}{Theorem}[section]

\newtheorem{definition}[theorem]{Definition}
\newtheorem{Lemma}[theorem]{Lemma}
\newtheorem{Proposition}[theorem]{Proposition}
\newtheorem{Remark}[theorem]{Remark}
\newtheorem{corollary}[theorem]{Corollary}



\vglue50pt

\centerline{\Large\bf Convergence in law for }

\bigskip

\centerline{\Large\bf the branching random walk seen from its tip}

\bigskip
\bigskip

\centerline{by}

\medskip

\centerline{Thomas Madaule}

\medskip

\centerline{\it Universit\'e Paris XIII}

\bigskip
\bigskip
\bigskip

{\leftskip=2truecm \rightskip=2truecm \baselineskip=15pt \small

\noindent{\slshape\bfseries Abstract.} Considering a critical branching random walk on the real line. In a recent paper, Aïdékon \cite{Aid11} developed a powerful method to obtain the convergence in law of its minimum after a log-factor translation. By an adaptation of this method, we show that the point process formed by the branching random walk seen from the minimum converges in law to a decorated Poisson  point process. This result,  confirming a  conjecture of Brunet and Derrida \cite{BDe11},  can be viewed as a discrete analog of the corresponding results for the branching brownian motion, previously established by  Arguin et al. \cite{ABK10} \cite{ABK11} and Aïdékon et al. \cite{ABBS11}.

} 

\bigskip
\bigskip

\section{Introduction}
We consider a branching random walk on the real line $\r$. Initially, a single particle sits at the origin. At time $1$, the particle gives birth to some children which form the first generation of the branching random walk and whose positions are given by a point process $L$ on $\r$. At time $2$, each of the particles in the first generation gives birth to new particles that are positioned --with respect to their birth positions-- according to the law of the same point process $L$; they form the second generation.  And so on. We assume that, each particle produces new particles independently of other particles in the same generation, and of everything up to that generation.

Let $\mathbb{T}$ be the genealogical tree of the particles in the branching random walk, then $\mathbb{T}$ is a Galton-Watson tree.  We write $|z|=n$ if a particle  $z$  is in the $n$-th generation, and denote its position by $V(z)$.\nomenclature{$\mathbb{T}$}{: the genealogical tree of particles} \nomenclature{$|z|$}{: generation of $z\in \mathbb{T}$} The collection of positions $(V(z),z\in \mathbb{T})$ is our branching random walk.\nomenclature{$V(z)$}{: position of $z\in \mathbb{T}$ }

The study of the minimal position $M_n:=\min_{\vert z \vert =n}    V(z)$   \nomenclature{$M_n$}{: minimal position of the branching random walk at time $n$}  has attracted many recent interests.  The law of large numbers for the speed of the minimum goes back to the works of Hammersley \cite{Ham74}, Kingman \cite{Kin75} and Biggins \cite{Big77}. The second order problem was recently studied separately by Hu and Shi \cite{HuS09} (a.s fluctuation), and Addario-Berry and Reed \cite{Add09}. In \cite{Add09}, the authors computed the expectation of $M_n$ up to  $O(1)$, and showed, under suitable assumptions, that the sequence of the minimum is tight around its mean. Through recursive equations, Bramson and Zeitouni \cite{BZe09} obtained the tightness of $M_n$ around its median, assuming some hypotheses on the decay of the tail distribution. A definitive answer was recently given by  Aïdékon \cite{Aid11}, where he  proved the convergence of the minimum $M_n$ centered around $\frac{3}{2}\log n$ for the general class of critical branching random walks.

One problem of great interest in the study of branching random walk is to characterize  its behaviour seen from the minimal position, namely, the asymptotic of the point process formed by $\{ V(z)-M_n, |z|=n\}$ as $n \to \infty$. The corresponding problem for the branching Brownian motion (the continuous analogue of branching random walk) was solved very recently by Arguin, Bovier, Kistler \cite{ABK10}, \cite{ABK11} and in parallel by Aïdékon, Beresticky, Brunet, Shi \cite{ABBS11}.

The aim of this paper is to establish the analogous results for branching random walk. Our main result, summed up in Theorem 1.1, will give the existence of the limiting point process together with a partial description, which also confirms the prediction in Brunet and Derrida \cite{BDe11}. Our method, largely inspired by Aïdékon \cite{Aid11}, consists in an analysis of the Laplace transform of the point process.

Following \cite{Aid11}, we assume\nomenclature{$\P$}{: law of the branching random walk}
\begin{equation}
\E\left[\underset{|z|=1}{\sum}1\right]>1,\qquad \E\left[\underset{|z|=1}{\sum}\ee^{-V(z)}\right]=1,\qquad \E\left[\underset{|z|=1}{\sum}V(z)\ee^{-V(z)}\right]=0.
\label{boundarycase}
\end{equation}
Every branching random walk satisfying mild assumptions can be reduced to this case by a linear transformation. We refer to the Appendix A in \cite{Jaf09} for a precise discussion. Notice that we allow $\E\Big[\underset{|z|=1}{\sum}1\Big]=\infty$, and even $\P\Big(\underset{|z|=1}{\sum} 1=\infty \Big)>0$. The couple $(M_n, W_{n,\beta})$ is the most often encountered random variables in our work, with \nomenclature{$W_{n,\beta}$}{: partition function of the branching random walk}
$$M_n:=\min\{V(x),|z|=n\},\qquad W_{n,\beta}:=\underset{|z|=n}{\sum}\ee^{-\beta V(z)},\qquad \beta>1,\,n\geq 1.$$

We also need the \nomenclature{$Z_n$}{$:=\underset{|z|=n}{\sum}V(z)\ee^{-V(z)}$ \nomenclature{$Z_\infty$}{$:=\underset{n\to\infty}{\lim}Z_n$}, the derivative martingale} \noindent{\it derivative martingale}
\begin{equation}
Z_n := \underset{|z|=n}{\sum} V(z)\ee^{-V(z)},\qquad Z_\infty=\underset{n\to\infty}{\lim} Z_n.
\end{equation}
By \cite{BKy04} and \cite{Aid11} we know that $Z_\infty$ exists almost surely and is strictly positive on the set of non extinction of $\mathbb{T}$. As in the continuous case \cite{ABBS11}, we introduce the point process formed by the particles of the recentered branching random walk: \nomenclature{$\mu_n$}{: point process seen from its tip} 
$$\mu_n:=\underset{|z|=n}{\sum}\delta_{\{V(z)-\frac{3}{2}\log n+\log Z_\infty\}}, \quad n\geq 1. $$

As in the setting of Kallenberg \cite{Kal76} (pp 10), a point process is considered as a random map with value in $\mathcal{N}$ the space of locally finite counting measure equipped with the vague topology. Then convergence in distribution of point processes will mean weak convergence of the corresponding distributions with respect to the vague topology (see \cite{Kal76} pp 42). 

We will show the existence in distribution of a limiting point process of $\mu_n$ as $n\to \infty$, from which we deduce results on \nomenclature{$\mu_n'$}{: point process seen from its minimal position} $\mu_n':=\underset{|z|=n}{\sum}\delta_{\{V(z)-M_n\}}$, $n\geq 1$. 

Writing for $y\in \r\cup\{\infty\}$, $y_+:= \max(y,0)$, we introduce the random variables \nomenclature{$X$, $\widetilde{X}$}{: defined in (\ref{intcdtn})}
\begin{equation}
X:=\underset{|z|=1}{\sum}\ee^{-V(z)},\qquad \widetilde{X}:=\underset{|z|=1}{\sum}V(z)_+\ee^{-V(z)}.
\label{intcdtn}
\end{equation}
with the convention $\infty\ee^{-\infty}=0$. We finally assume that 
\begin{eqnarray}
\label{1.3bis}
&&\text{the distribution of }L\text{ is non-lattice},
\\
&&\E\Big[\underset{|z|=1}{\sum}V(z)^2\ee^{-V(z)}\Big]<\infty,\quad\E\Big[X(\log_+ X)^2\Big]<\infty,\quad \E\Big[\tilde{X}(\log_+ \tilde{X})\Big]<\infty.
\label{condint}
\end{eqnarray}

The main result of this paper is the following theorem: 
\begin{theorem}
Under (\ref{boundarycase}), (\ref{1.3bis}) and (\ref{condint}), as $n\to \infty$, conditioning on the set of non-extinction, the pair $(\mu_n,Z_n)$ converges jointly in distribution to $(\mu_\infty,Z_\infty)$ where $\mu_\infty$ and $Z_\infty$ are independent and $\mu_\infty$ is obtained as follows.\nomenclature{$\mu_\infty$}{: limiting law of $\mu_n$}

(i) Define $\mathcal{P}$ to be \nomenclature{$\mathcal{P}$}{: decorated Poisson point process} a Poisson point process on $\r$, with intensity measure $\lambda\ee^x dx$ for some positive real constant $\lambda$.

(ii) For each atom x of $\mathcal{P}$, we attach a point process $x+\mathcal{D}^{(x)}$ where $\mathcal{D}^{(x)}$ are independent copies of a certain point process $\mathcal{D}$ in $\r_+$.\nomenclature{$D$, $D^{(x)}$}{: decoration of $\mathcal{P}$}

(iii) The point process $\mu_\infty$ is the superposition of all the point processes $x+\mathcal{D}^{(x)}$, i.e, $\mu_\infty:=\{x+y:x\in \mathcal{P},y\in \mathcal{D}^{(x)}\}$.
\label{maintheorem}
\end{theorem}
\paragraph{Remark:} The point process $\mu_\infty$ is called "decorated Poisson point process with decoration $\mathcal{D}$". We refer to \cite{Mai11} for a more complete description. 
\begin{corollary}
Under (\ref{boundarycase}), (\ref{1.3bis}) and (\ref{condint}), conditioning on the set of non-extinction, seen from the leftmost particle, the point process $\mu_n'$ formed by the particles $\{V(u)-M_n, \vert u\vert =n \}$ converges in distribution to the point process $\mu_\infty'$ obtained by replacing the Poisson point process $\mathcal{P}$ in step (i) above by $\mathcal{P}'$ described in step (i)' below:

(i)' Let ${\bm \ee}$ be a standard exponential random variable. Conditionally on ${\bm \ee}$, define $\mathcal{P}'$ to be a Poisson point process on $\r_+$, with intensity measure ${\bm \ee}\, \ee^x\1_{\r_+}dx$ to which we add an atom in $0$.

The decoration point process $\mathcal{D}$ remains the same.
\end{corollary}

These two results imitate the corresponding results for the branching Brownian motion, in particular Theorem 2.1 and Corollary 2.2 of  Aïdékon, Beresticky, Brunet and Shi \cite{ABBS11} (and also those of \cite{ABK10} and \cite{ABK11}).  However, we do not adopt the same method as in \cite{ABBS11}  because, firstly the spine decomposition for the branching random walk leads to a use of Palm measures, which is much more complicated than in the case of branching Brownian motion, and secondly, the path decomposition for a random walk is also more complex than in the Brownian case. Instead, we shall imitate the fine analysis of Aïdékon \cite{Aid11} to study the Laplace transform of $\mu_n$.  
More precisely,  the main step in the proof of Theorem \ref{maintheorem}  is to establish the convergence in law of $(n^{\frac{3}{2}\beta_1}W_{n,\beta_1},...,n^{\frac{3}{2}\beta_k}W_{n,\beta_k})$ for any $k\geq 1$ and any  $\beta_k>...>\beta_1>1$. A crucial observation, inspired from \cite{Aid11}, is that the convergence in law of this $k$-dimensional random vector can be reduced to the study of its tail distribution.  From this, we can prove the convergence in law stated in Theorem  \ref{maintheorem}, and  as a by-product, we shall also get an expression for the Laplace transform of the limiting point process. The later  might have some independent interest for further analysis of $\mu_\infty$. 

Note that the present paper provides a simplification of the method of \cite{Aid11} and could be used to recover the Theorem 1.1 of \cite{Aid11}. Indeed using the deep understanding induced by the work of \cite{Aid11} we are able to skip the use of the ''killed branching random walk" (see \cite{Aid11} for a definition). However this simplification necessitates to prove numerous lemmas which have very close analogue in \cite{Aid11}. Most often these lemmas are stated and proved in Appendices A, B and C.

The paper is organized as follows: Section 2 contains the main estimates on the tail of distribution of $(n^{\frac{3}{2}\beta_1}W_{n,\beta_1},...,n^{\frac{3}{2}\beta_k}W_{n,\beta_k})$ for $k\geq 1$ and $\beta_k>...>\beta_1>1$, from which we establish the convergence of some Laplace transforms of $\mu_n$ (Theorem \ref{laplace}) and give the proof of  Theorem \ref{maintheorem}.  Section 3 is devoted to the proof Theorem \ref{laplace} by admitting Proposition \ref{Propconclusion2}. Finally, we prove in Sections 4 and 5 respectively  Propositions \ref{lemdomination} and \ref{Propconclusion2}.

\section{Main steps of the proof of Theorem \ref{maintheorem}}
To shorten the statements we introduce some notations: For $n\geq 1$, $\beta>1$, define \nomenclature{$\widetilde{M}_n,\, \widetilde{W}_{n,\beta},\,\widetilde{W}_{n-|u|,\beta}$}{: respectively $M_n-\frac{3}{2}\log n$ and  $n^{\frac{3}{2}\beta}\times$ $W_{n,\beta},\,W_{n-|u|,\beta}$}
\begin{eqnarray*}
&\widetilde{W}_{n,\beta}:=n^{\frac{3}{2}\beta}W_{n,\beta},\qquad \widehat{\mu}_n(\beta)=n^{\frac{3}{2}\beta}\underset{|z|=n}{\sum}\ee^{-\beta (V(z)+\log Z_\infty)}.&
\end{eqnarray*}
Remark that $\widehat{\mu}_n(\beta)$ is also equal to $\int_\r \ee^{-\beta x}d\mu_n(x)$. In a general context many quantities   with  tilde are associated with the natural  normalization $n^{\frac{3}{2}\beta}$ except for some obvious abuse of notation: For example in the sequel we will denote for brevity  $\widetilde{W}_{n-|u|,\beta}:= n^{\frac{3}{2}\beta}W_{n-|u|,\beta}$. In a similar spirit we write $\widetilde{M}_n:=M_n-\frac{3}{2}\log n $ and $\widetilde{M}_{n-|u|}:= M_{n-|u|}-\frac{3}{2}\log n$ for some vertex $\vert u \vert   \le n$ (we shall recall these notations to avoid any confusion). Throughout the paper $\N^*$, $\r_+$ and $\r^{*}_+$ denote respectively $\N_{\backslash 0}$, $[0,\infty)$ and $(0,\infty)$.  At last we often encounter notations ${\bm \delta}$, ${\bm \beta}$ and ${\bm y}$ for respectively $(\delta_1,...,\delta_k)$, $(\beta_1,...,\beta_k)$ and $(y_1,...,y_k)$, with some $k\geq 1$ determined in the context.

\subsection{Main preliminary results}
In this section we state the main results which will lead to the proof of Theorem \ref{maintheorem} (deferring their proofs to the next sections).

The following Proposition give a uniform control on the tail of distribution of the process $(n^\frac{3\beta}{2}W_{n,\beta})_{n\geq 1}$.
\begin{Proposition}
\label{lemdomination}
Under (\ref{boundarycase}) and (\ref{condint}), there exists $c_{1}>0$ such that for any $n>1$, and $x\geq 1$,
\begin{equation}
\label{eqlemdomination}
\P\left(\widetilde{W}_{n,\beta}\geq \ee^{\beta x}\right)\leq c_{1} (1+x)\ee^{-x}.
\end{equation}
\end{Proposition}

\begin{Proposition}
\label{Propconclusion2}
Under (\ref{boundarycase}) and (\ref{condint}), for any $d\in \N^* $ and ${\bm \beta}\in (1,\infty)^d$ there exists a function $\rho_{\bm \beta}:{\bm \theta}\ni(\r_{+}^*)^{d}\to \rho_{\bm \beta}({\bm \theta})\in (0,\infty)$, which satisfies the following:\nomenclature{$ \rho({\bm \theta})$}{: defined in (\ref{queuededistrib})}
For any ${\bm \theta}\in \r_{+,*}^d,\,  \epsilon>0$,  there exists  $(A,N)_{(\epsilon)}\in \r_+\times \N$  such that $\forall  n>N$ and $ x\in[ A,\frac{3}{2}\log n-A]$, we have
\begin{equation}
\label{sansDelta2}
\left|\frac{\ee^x}{x}\E\left(1-\exp\{-\sum_{i=1}^d \theta_i \ee^{-\beta_i x} \widetilde{W}_{n,\beta_i}\} \right)-\rho_{\bm \beta}({\bm \theta })\right|\leq \epsilon. 
\end{equation}
Moreover $\lim_{{\bm \theta }\to 0}\rho_{\bm \beta}({\bm \theta})=0$ and for any $y\in \r,\, {\bm \theta} \in (\r_{+}^*)^d$, $\rho_{\bm \beta}(\theta_1\ee^{\beta_1 y},..., \theta_d\ee^{\beta_d y})= \ee^y\rho_{\bm \beta}({\bm \theta})$ (this last equality is plain by a change of variable in (\ref{sansDelta2})).
\end{Proposition}

The key step in the proof of Theorem \ref{maintheorem} is the following result: 
\begin{theorem}
\label{laplace}
Under (\ref{boundarycase}) and (\ref{condint}), $\forall d\in \N,\, {\bm \beta}\in(1,\infty)^d$,

\begin{equation}
\label{eqlaplace}
\forall  \alpha\in \r_+,\,\,\underset{n\to \infty}{\lim}\E\left(\ee^{-\underset{i=1}{\overset{d}{\sum}}\theta_i \widehat{\mu}_n(\beta_i)}\ee^{-\alpha Z_\infty}\1_{\{Z_\infty>0\}}\right)=\ee^{-\rho_{\bm \beta}({\bm \theta}) }\E\left(\ee^{-\alpha Z_\infty}\1_{\{Z_\infty>0\}}\right).
\end{equation}
In particular as $n\to \infty$, conditionally on $\{Z_\infty >0\}$,  $(\widehat{\mu}_n(\beta_1),...,\widehat{\mu}_n(\beta_d))$ converges in law to some random vector $(\widehat{\mu}_\infty(\beta_1),...,\widehat{\mu}_\infty(\beta_d))$ independent of $Z_\infty$.
\end{theorem}

\subsection{Proof of the Theorem \ref{maintheorem} by admitting Theorem 2.3}
Let us introduce the conditional probability $\P^{*}(.):=\P(.|\text{non-extinction})$. Recall that under $\P^{*}$, $Z_\infty>0$ a.s. To prove the Theorem \ref{maintheorem} we have to keep in mind two facts:

-According to Theorem \ref{laplace}, for any $l\in \N^*$ and $\bbeta\in (\{2,3,...\})^l$ the vector $(\widehat{\mu}_n(\beta_1),...,\widehat{\mu}_n(\beta_l) )$ converges in law under $\P^{*}$. We deduce that for any $\ttheta\in \r^l $, $\overset{l}{\underset{i=1}{\sum}}\theta_i\widehat{\mu}_n(\beta_i)$ converges also in law. But $\overset{l}{\underset{i=1}{\sum}}\theta_i\widehat{\mu}_n(\beta_i)$ is the same as $\int_\r Q(\ee^{-x})\d\mu_n(x)$ with $Q(X):=\overset{l}{\underset{i=1}{\sum}}\theta_i X^{\beta_i}$. Hence, if $Q$ is polynomial function such that $Q(0)=Q'(0)=0$, then $\int_\r Q(\ee^{-x})\d\mu_n(x)$ converges in law in law under $\P^{*}$.

-In \cite{Aid11}, Elie Aïdékon has proved that under $\P^{*}$, $M_n-\frac{3}{2}\log n$ converges in law but $M_n-\frac{3}{2}\log n + Z_\infty$ is nothing but the smallest atom of the point process $\mu_n$. Thus it is clear that $\underset{b\to \infty}{\lim}\underset{n\in\N}{\sup}\,\P^{*}\left(\mu_n(-\infty,-b]>0\right)=0$.

Let $C_c(\r)$ be the set of continuous functions in $\r$ with compact support.\nomenclature{$C_c(\r)$}{: set of continuous functions with compact support} The existence of $\mu_\infty$ is now a consequence of Kallenberg \cite{Kal76}. Indeed the Lemma 5.1 in \cite{Kal76} says that $\mu_n$ converges in law to some $\mu_\infty$ (for the vague topology see \cite{Kal76} chapter 4) providing that $\forall f\in C_c(\r)$, $\left( \int_\r f(x)d\mu_n(x) \right)_{n\in \N}$ converges in law to some random variable $\mu(f)$. 
\\

Without loss of generality suppose that $f(x) =g(x) \ee^{-2 x}$ with $g\in C_c(\r)$. Let $\epsilon>0$. Let $b\in \r$ such that $g(x)=0$ for any $x\in [-b,b]^c$ and $\underset{n\in\N}{\sup}\,\P^{*}\left(\mu_n(-\infty,-b]>0\right)\leq \epsilon$. According to Stone-Weierstrass' Theorem there exists a sequence of polynomial function $Q_q\in \r[x]$ such that $\underset{y\in(0,\ee^{b}]}{\sup}\left|Q_q(y)-g(\log \frac{1}{y})\right|\leq \frac{1}{q}$. By a change of variable this is equivalent to $ \underset{y\in[-b,+\infty)}{\sup}\left|Q_q(\ee^{-y})-g(y)\right|\leq \frac{1}{q} $.

For $\theta\in \r$ and $n,\,p,\, q\in \N^*$, the triangle inequality implies that

\begin{eqnarray}
\label{inepousser}
\left|\E^{*}\left(\ee^{i\theta \int_\r g(x)\ee^{-2x}d\mu_n(x)}\right)-\E^{*}\left(\ee^{i\theta \int_\r g(x)\ee^{-2x}d\mu_p(x)}\right)\right|\leq (1)_{n,q}+(1)_{p,q}+(2)_{n,p,q},
\end{eqnarray}
with
\begin{eqnarray*}
&&(1)_{n,q}:=\left|1-\E^{*}\left(\ee^{i\theta \int_\r [g(x)-Q_q(\ee^{-x})]\ee^{-2x}d\mu_n(x)}\right)\right|,
\\
&&(2)_{n,p,q}:=\left|\E^{*}\left(\ee^{i\theta \int_\r Q_q(x)\ee^{-2x}d\mu_n(x)}\right)-\E^{*}\left(\ee^{i\theta \int_\r Q_q(x)\ee^{-2x}d\mu_p(x)}\right)\right|.
\end{eqnarray*}
For any $y\in \r$, $|1-\ee^{iy}|\leq |y|$, therefore for any $n$ and $q\in \N^*$,

\begin{eqnarray*}
(1)_{n,q}&\leq& 2\P^{*}\left(\int_\r\ee^{-2x}d\mu_n(x)\geq \sqrt{q}\right)+2\P^{*}(\mu_n(-\infty,-b]>0)+
\\
\qquad \qquad&& \E^{*}\left(\left(\ee^{i\theta \int_\r [g(x)-Q_q(\ee^{-x})]\ee^{-2x}d\mu_n(x)}-1\right)\1_{\{\int_\r \ee^{-2x}d\mu_n(x)\leq \sqrt{q},\,\mu_n(-\infty,-b]=0 \}}\right)
\\
&\leq& 2\P^{*}\left(\widehat{\mu}_n(2)\geq \sqrt{q}\right)+2\epsilon+|\theta|\frac{\sqrt{q}}{q}.
\end{eqnarray*}
Thanks to the tightness of $( \widehat{\mu}_n(2))_{n\in \N}$, we can choose $q_0$ sufficiently large such that $\underset{n\in\N}{\sup}\,(1)_{n,q_0}\leq 3\epsilon$. As $\int_\r Q_{q_0}(\ee^{-x})d\mu_n(x)$ converges in law, for $n$ and $p$ sufficiently large $(2)_{n,p,q_0}\leq \epsilon$ uniformly in $\theta$ in any compact set. Thus the sequence $\E^{*}\left(\ee^{i\theta \int_\r f(x)\ee^{-2x}d\mu_n(x)}\right)$ satisfies Cauchy's criterion and hence admits a limit that we denote $\Psi_f(\theta)$. Moreover the convergence is uniform on every compact in $\theta$, thus $\Psi_f(\theta)$ is continuous at $0$.
\\

We have proved that the limiting law $\mu_\infty$ exists. To obtain the description of this point process as a decorated Poisson point process it suffices, according to Corollary 5.2 of  Maillard \cite{Mai11} (see also \cite{DMZ05}), to prove that $\mu_\infty$ is {\it superposable}. See \cite{BDe11} for the origin of this idea. We recall what this notion means:

Let $\mathcal{N}$ be the space of locally finite counting measures on $\r$.\nomenclature{$ \mathcal{N}$}{: space of bounded counting measures} For every $x\in \r$ define the translation operator $T_x:\mathcal{N}\to\mathcal{N}$, by $(T_x\mu)(A)=\mu(A+x)$ for every Borel set $A\subset \r$.\nomenclature{$ T_x$}{: translation operator} Let $\mathcal{L}'$ be an independent copy on $\mathcal{L}$. We say that $\mathcal{L}$ is $superposable$, if

$$T_\alpha \mathcal{L}+T_\beta \mathcal{L}'\overset{(d)}{=}\mathcal{L},\text{ for every }\alpha,\beta\in\r \text{ such that }\ee^{-\alpha}+\ee^{-\beta}=1.$$
In view to prove the {\it superposability} of $\mu_\infty$, for any $(a,b)\in \r^2$, we introduce $\mathbb{T}^{a}$ and $\mathbb{T}^{b}$ the genealogical trees formed by two independent branching random walks starting respectively at $a$ and $b$. Following section 1 we introduce all the objects related to $\mathbb{T}^{a}$ and $\mathbb{T}^b$ by adding an extra superscript $a$ or $b$ in the notation (i.e $Z_n^a:= \sum_{|z|=n,\, z\in \mathbb{T}^a}V(z)\ee^{-V(z)}$ , $Z_n^b:= \sum_{|z|=n,\, z\in \mathbb{T}^b}V(z)\ee^{-V(z)}, \mu_n^a, \mu_n^b...$). We also define
\begin{eqnarray*}
&&\mu_n^{a,b}:=\mu_n^a+ \mu_n^b=\underset{u\in \mathbb{T}^{a},\, \vert u\vert =n}{\sum}\delta_{V(u)+Z_\infty^{a}}+\underset{u\in \mathbb{T}^{b},\, \vert u\vert =n}{\sum}\delta_{V(u)+Z_\infty^{b}},\text{   and }
\\
&& \widehat{\mu}_n^{a, b}(\beta):=\int_\r \ee^{-\beta x} d\mu_n^{a,b}(x),\quad \beta>1.
\end{eqnarray*}
We note that $\widehat{\mu}_n^{a, b}(\beta)\overset{(d)}{=}\ee^{-\beta a}\widehat{\mu}_n^{(1)}(\beta)+\ee^{-\beta b}\widehat{\mu}_n^{(2)}(\beta)$ with $\widehat{\mu}_n^{(i)}(\beta)$, $i\in \{1,2\}$, two independent copy of $\widehat{\mu}_n(\beta)$. The key point to prove that $\mu_\infty$ is {\it superposable} is the following Corollary:
\begin{corollary}
\label{superposabilite}
Under (\ref{boundarycase}) and (\ref{condint}), for any $d\in \N^*,\,{\bm \beta}\in(1,\infty)^d$, $a$ and $b\in\r$ such that $\ee^{-a}+\ee^{-b}=1$ we have
\begin{equation}
\label{eqlaplace2}
\underset{n\to \infty}{\lim}\E\left(\ee^{ -\underset{i=1}{\overset{d}{\sum}}\theta_i \widehat{\mu}_n^{a,b}(\beta_i) }\1_{\{Z^a_\infty>0,\,Z^b_\infty>0\}}\right)=\ee^{-\rho_{\bm \beta}({\bm \theta})}\P(Z^a_\infty>0,\,Z^b_\infty>0),
\end{equation}
where $\rho_{\bm \beta}$ is the same function in (\ref{eqlaplace}) and (\ref{eqlaplace2}). Therefore when $n\to \infty$ the limit in law of  $(\widehat{\mu}^{a,b}_n(\beta_1),...,\widehat{\mu}^{a,b}_n(\beta_d))$ conditionally on $\{Z^a_\infty>0,\,Z^b_\infty>0\}$ is also the law of  $(\widehat{\mu}_\infty(\beta_1),...,\widehat{\mu}_\infty(\beta_d))$ (see Theorem \ref{laplace}). 
\end{corollary}
\noindent {\it Proof of Corollary \ref{superposabilite}.} We apply Theorem \ref{laplace} to obtain $\underset{n\to \infty}{\lim} \E\Big(\ee^{ -\underset{i=1}{\overset{d}{\sum}}\theta_i \widehat{\mu}_n^{a,b}(\beta_i) }  \1_{\{Z^a_\infty>0,\,Z^b_\infty>0\}}\Big) = \ee^{-\rho^a_{\bm \beta}({\bm \theta})-\rho^b_{\bm \beta}({\bm \theta})} \P(Z^a_\infty>0,\,Z^b_\infty>0)$ with $\forall x\in \r$, $\rho^x_{\bm \beta}({\bm \theta}):=\rho_{\bm \beta}(\ee^{-\beta_1 x}\theta_1,...,\ee^{-\beta_l x}\theta_d)$. By the Proposition \ref{Propconclusion2}, $\rho^x_{\bm \beta}({\bm \theta})= \ee^{-x}\rho_{\bm \beta}({\bm \theta})$. Applied to $x=a$ and $x=b$ this gives the result.
\hfill$\Box$
\\

\noindent{\it Proof of the superposability of $\mu_\infty$.} Let $(a,b)\in\r^2$ be two constants such that $\ee^{-a}+\ee^{-b}=1$, via the proof of the existence of $\mu_\infty$ and Corollary \ref{superposabilite}, it is clear that conditionally on $\{Z^a_\infty>0,\,Z^b_\infty>0\}$, $\mu_n^{a,b}$ converges also in law to $\mu_\infty$. Moreover for any $n\in \N$,  $\mu_n^{a,b}= \mu_n^a+\mu_n^b\overset{\text{law}}{=}T_a\mu_n^{(1)}+T_b\mu_n^{(2)}$ with $\mu_n^{(i)}$, $i\in \{1,2\}$ two independent branching random walks. So if we reformulate in terms of {\it superposability} we have
\begin{eqnarray*}
T_a\mu_\infty^{(1)}+T_b\mu_\infty^{(2)}\overset{(d)}{=}\underset{n\to\infty}{\lim}\mu_n^{a,b}\overset{(\text{law})}{=}\mu_\infty.
\end{eqnarray*}
thus $\mu_\infty$ is a point process {\it superposable}.
\hfill $\Box$

Assuming the Theorem \ref{laplace} the proof of Theorem \ref{maintheorem} is complete.\hfill$\Box$

\section{Proof of Theorem \ref{laplace} by admitting Proposition \ref{Propconclusion2}} 
When $Y$ is a non negative random variable and $\Xi$ an event, we often will write $\E(Y;\Xi)$ for $\E(Y \1_{\Xi})$.

For  any vertex $z\in \mathbb{T}$ we denote by $[\emptyset,z]$ the unique shortest path relating $z$ to the root $\emptyset$, and $z_i$ (for $i<|z|$) the vertex on $[\emptyset,z]$ such that $|z_i|=i$. The trajectory of $z\in \mathbb{T}$, $|z|=n$, corresponds to the ancestor's positions of $z$, i.e the vector $(V(z_1),...,V(z_n))$. If $x$ is an ancestor of $y$ we will write $x<y$.

Let $d\in \N^*$, $({\bm \beta},{\bm \theta},\alpha)\in (1,+\infty)^d\times (\r_{+}^*)^{d}\times \r_{+}^{*}$.

Let $\mathcal{Z}[A]$ (and $\mathcal{F}_A$ the corresponding $\sigma$-field) be the set of particles absorbed at level $A\in \r_+$, i.e.\nomenclature{$ \mathcal{Z}[A]$}{: set of particles}\nomenclature{$\mathcal{F}_A$}{: sigma field generated by the particles in $\mathcal{Z}[A]$} 
$$\mathcal{Z}[A]:=\{u\in \mathbb{T}: V(u)\geq A, V(u_k)<A \text{ }\forall k<|u|\}\text{ and } \mathcal{F}_A:=\sigma\left( (u,V(u));\,\,\,u\in \mathcal{Z}[A]\right).$$
With a slight extension of \cite{BKy04}, equality (5.2) pp 46 in \cite{Aid11} affirms that $Z_A:=\underset{u\in \mathcal{Z}[A]}{\sum}V(u)\ee^{-V(u)}$ satisfies that\nomenclature{$Z_A$}{: defined in (\ref{Z[A]})}
\begin{equation}
\label{Z[A]}
\underset{A\to \infty}{\lim} Z_A=Z_\infty  \qquad a.s.
\end{equation}
In a first step we will show that
\begin{Lemma}
Under (\ref{boundarycase}) and (\ref{condint}),
\label{propfi}
\begin{equation}
\label{finnn}
\underset{A\to \infty}{\lim}\underset{n\to\infty}{\lim}\,\E\left(\ee^{-\overset{d}{\underset{i=1}{\sum}}\theta_i \ee^{-\beta_i \log Z_A} \widetilde{W}_{n,\beta_i}}\1_{\{Z_A>0\}}\ee^{-\alpha Z_A}\right)=\ee^{-\rho_{\bm \beta}({\bm \theta})}\E(\ee^{-\alpha Z_\infty};\,Z_\infty>0).
\end{equation}
\end{Lemma}
{\noindent \it Proof of Lemma \ref{propfi}.}
For every $n\in\N$ and $L>0$, we define $\Xi_{A}(n,L)\equiv\Xi_A$ by\nomenclature{$\Xi_{A}$}{: set of particles in the branching random walk}
\begin{equation}
\Xi_{A}(n,L):=\{\underset{u\in \mathcal{Z}[A]}{\max}|u|\leq (\log n)^{10},\underset{u\in \mathcal{Z}[A]}{\max} V(u)\leq A + \frac{1}{30}\log n,\, \log Z_A\in [-L,L]\}.
\label{defchiii}
\end{equation}
On the set of non-extinction $Z_\infty>0$ and $ M_n {\to} +\infty$ a.s, thus, conditional on $Z_\infty>0$, the probability of $\Xi_A$ increases to $1$ when $n$, then $A$ and then $L$ go to infinity. 

At first we study \nomenclature{$\mathcal{F}_\Xi$}{: defined in (\ref{defFchii})}
\begin{equation}
\mathfrak{F}_\Xi(A,n,L)=\mathfrak{F}_\Xi:= \E\Big(\ee^{-\overset{d}{\underset{i=1}{\sum}}\theta_i \ee^{-\beta_i \log Z_A}\widetilde{W}_{n,\beta_i}}\ee^{-\alpha Z_A}\1_{\{Z_A>0\}};\,\Xi_A\Big).
\label{defFchii}
\end{equation}
On $\Xi_A$ we have $\widetilde{W}_{n,\beta}:=\underset{u\in \mathcal{Z}[A]}{\sum}\ee^{-\beta V(u)}\widetilde{W}_{n,\beta}^u$ with $W_{n,\beta}^u:=\underset{z>u,\vert z\vert =n}{\sum}\ee^{-\beta (V(z)-V(u))}$ (recalling that $\widetilde{W}_{n,\beta}^u=n^{\frac{3}{2}\beta}W_{n,\beta}^u$). For any $\beta>1$, $W_{n,\beta}^u$ has the same law has $W_{n-|u|,\beta}$, then the Markov property leads to 
\begin{eqnarray*}
\mathfrak{F}_\Xi&=&\E\Big(\ee^{-\overset{d}{\underset{i=1}{\sum}}\theta_i \underset{u\in \mathcal{Z}[A]}{\sum}\ee^{-\beta_i (V(u)+\log Z_A)}\widetilde{W}^u_{n,\beta_i}}\ee^{-\alpha Z_A}\1_{\{Z_A>0\}};\Xi_A\Big)
\\
&=&\E\Big(\underset{u\in \mathcal{Z}[A]}{\prod}\E\Big(\ee^{-\overset{d}{\underset{i=1}{\sum}}\theta_i \ee^{-\beta_i (V(u)+\log Z_A)}n^{\frac{3\beta}{2}}W_{n-|u|,\beta_i}}\Big|\mathcal{F}_A\Big)\ee^{-\alpha Z_A}\1_{\{Z_A>0\}} ; \Xi_A\Big)
\\
&=&\E\Big(\underset{u\in \mathcal{Z}[A]}{\prod} \E\big(\ee^{-\overset{d}{\underset{i=1}{\sum}}\theta_i \ee^{-\beta_i (V(u)+\log Z_A)}n^{\frac{3\beta}{2}}W_{n-|u|,\beta_i}}\big)  \ee^{-\alpha Z_A}\1_{\{Z_A>0\}} ; \Xi_A\Big).
\end{eqnarray*}
On $\Xi_A$, for any $u\in \mathcal{Z}[A]$,  $V(u)+\log Z_A\in [A-L,L+A+\frac{1}{30}\log n]$, then by the Proposition \ref{Propconclusion2}, there exists $A,N$ large enough such that $\forall n>N$, we have for any $u\in \mathcal{Z}[A]$, 
\begin{equation}
\label{prebecom} \Big|1-\E\big(\ee^{-\overset{d}{\underset{i=1}{\sum}}\theta_i \ee^{-\beta_i( V(u)+\log Z_A)}n^{\frac{3\beta}{2}}W_{n-|u|,\beta_i}} \Big|\mathcal{F}_A \big)-\rho_{\bm \beta}( {\bm \theta})(V(u)+\log Z_A)\frac{\ee^{-V(u)}}{Z_A}\Big|\leq \epsilon  (V(u)+\log Z_A) \frac{\ee^{-V(u)}}{Z_A}.
\end{equation}
Moreover, with $A$ large enough such that $\rho_{\bf \beta}({\bf \theta})\frac{\log Z_A}{V(u)}\leq \rho_{\bf \beta}({\bf \theta})\frac{L}{A}\leq \epsilon$, (\ref{prebecom}) becomes
\begin{equation}
\label{becom} \Big|1-\E\big(\ee^{-\overset{d}{\underset{i=1}{\sum}}\theta_i \ee^{-\beta_i( V(u)+\log Z_A)}n^{\frac{3\beta}{2}}W_{n-|u|,\beta_i}}  \Big|\mathcal{F}_A\big)-\rho_{\bm \beta}(  {\bm \theta} )V(u)\frac{\ee^{-V(u)}}{Z_A}\Big|\leq 3\epsilon  V(u) \frac{\ee^{-V(u)}}{Z_A}.
\end{equation}
We deduce that
\begin{eqnarray*}
\mathfrak{F}_{\Xi}&\leq & \E\Big(\underset{u\in \mathcal{Z}[A]}{\prod}(1-[ \rho_{\bm \beta}({\bm \theta} )-3\epsilon ]\frac{V(u)}{Z_A}\ee^{-V(u)})\ee^{-\alpha Z_A}\1_{\{Z_A>0\}}\Big)
\\
&= &\E\Big(\ee^{\underset{u\in \mathcal{Z}[A]}{\sum}\log[1-[ \rho_{\bm \beta}({\bm \theta})-3\epsilon ]\frac{V(u)}{Z_A}\ee^{-V(u)}]}\ee^{-\alpha Z_A}\1_{\{Z_A>0\}}\Big)
\\
&\leq&\E\Big(\ee^{-(1-\epsilon)( \rho_{\bm \beta}( {\bm \theta})-3\epsilon)\frac{1}{Z_A}\underset{u\in Z[A]}{\sum}V(u)\ee^{-V(u)} }\ee^{-\alpha Z_A}\1_{\{Z_A>0\}}\Big)
\\
&=&\ee^{-(1-\epsilon)( \rho_{\bm \beta}( {\bm \theta})-3\epsilon) }\E\big(\ee^{-\alpha Z_A}\1_{\{Z_A>0\}}\big).
\end{eqnarray*}
Since $\underset{A\to\infty}{\lim}Z_A\overset{\text{a.s}}{=}Z_\infty$, we get 
$$\underset{A\to \infty}{\lim}\underset{n\to\infty}{\lim}\,\E\left(\ee^{-\overset{d}{\underset{i=1}{\sum}}\theta_i \ee^{-\beta_i \log Z_A} \widetilde{W}_{n,\beta_i}}\1_{\{Z_A>0\}}\ee^{-\alpha Z_A}\right)\leq \ee^{-\rho_{\bm \beta}({\bm \theta})}\E(\ee^{-\alpha Z_\infty};\,Z_\infty>0).$$ The lower bound follows similarly. This completes the proof of Lemma \ref{propfi}.
\hfill $\Box$

\medskip
\noindent{\it Proof of Theorem \ref{laplace}.} Because of (\ref{finnn}) it suffices to show that 
$$\E\Big(\ee^{-\overset{d}{\underset{i=1}{\sum}}\theta_i \ee^{-\beta_i \log Z_A}\widetilde{W}_{n,\beta_i}}\ee^{-\alpha Z_A}\1_{\{Z_A>0\}}\Big)\underset{A\to\infty}{\to}\E\Big(\ee^{-\overset{d}{\underset{i=1}{\sum}}\theta_i \widehat{\mu}_n(\beta_i)}\ee^{-\alpha Z_\infty}\1_{\{Z_\infty>0\}}\Big),$$
 uniformly on $n\in \N$. Let $\epsilon>0$. As $Z_A\overset{\text{p.s}}{\to}Z_\infty $ and the tightness of the sequence $(\widetilde{W}_{n,\beta})_{n\in \N}, \beta >1$ (see Proposition \ref{lemdomination}), there exists $K>0$ such that for any $n\geq 1$, 
 $$\P\big( \max_{i\in [1,d]}\theta_i \widetilde{W}_{n,\beta_i}\geq K \big)+ \P\big(\max_{i\in [1,d]}|\frac{1}{Z_A^{\beta_i}}-\frac{1}{Z_\infty^{\beta_i}} | + |Z_A-Z_\infty|\geq K,\, Z_\infty>0 \big)\leq \frac{\epsilon}{2}.$$
Then by employing the inequality $|\ee^{-x}-\ee^{-y}|\leq |x-y|$ for any $x,y \in \r_+$, uniformly in $n\geq 1$ we have
\begin{eqnarray*}
\big|\E\Big(\ee^{-\overset{d}{\underset{i=1}{\sum}}\theta_i \ee^{-\beta_i \log Z_A}\widetilde{W}_{n,\beta_i}}\ee^{-\alpha Z_A}\1_{\{Z_A>0\}}- \ee^{-\overset{d}{\underset{i=1}{\sum}}\theta_i \widehat{\mu}_n(\beta_i)}\ee^{-\alpha Z_\infty}\1_{\{Z_\infty>0\}}\Big)  \Big|\leq  |\P(Z_A>0)-\P(Z_\infty>0)|
\\
 +\frac{\epsilon}{2} + \E\big( K|\sum_{i=1}^d|\frac{1}{Z_A^{\beta_i}}-\frac{1}{Z_\infty^{\beta_i}} |  \1_{\{|\frac{1}{Z_A^{\beta_i}}-\frac{1}{Z_\infty^{\beta_i}} |\leq K \}} \big) + \E\big( \alpha |Z_A-Z_\infty |\1_{\{ |Z_A-Z_\infty|\leq K \}}\big).
\end{eqnarray*}
By dominated convergence this amount converges to $0$ when $A$ goes to infinity.

\hfill$\Box$

\section{Estimation on the tail of distribution of $\widetilde{W}_{n,\beta}$}

\subsection{The many-to-one formula  and Lyons' change of measure}
For $a\in \r$, we denote by $\P_a$ the probability distribution associated to the branching random walk starting from $a$, and $\E_a$ the corresponding expectation. Under (\ref{boundarycase}), we can define a random variable $X$ such that for any non-negative function $f$,
\begin{equation}
\E(f(X))=\E\left(\underset{|z|=1}{\sum}\ee^{-V(z)}f(V(z))\right).
\end{equation}
Moreover, via (\ref{condint}) we have $\sigma^2:=\E[X^2]<+\infty$. Let $(X_i)_{i\in \N^*}$ be a i.i.d sequence of copies of $X$. Write for any $n\in \N$, $S_n:=\underset{0<i\leq n}{\sum}X_i$, then $S$ is a mean-zero random-walk starting from the origin.\nomenclature{$(S_n)_{n\in\N}$}{: centred random walk associated to the many-to-one formula defined in (\ref{manytoone}) }\nomenclature{$\sigma^2 $}{$:=\E(S_1^2)$}
\begin{Lemma}[Biggins-Kyprianou]
Under (\ref{boundarycase}), for any $n\geq1$ and any measurable function $g:\r^n\to[0,+\infty)$,
\begin{equation}
\label{manytoone}
\E\left(\underset{|z|=n}{\sum}g(V(z_1),...,V(z_n))\right)=\E\left(\ee^{S_n}g(S_1,...,S_n)\right).
\end{equation}
\end{Lemma}
We can see the so-called many-to-one formula (\ref{manytoone}) as a consequence of Proposition \ref{lyons} below. Let introduce the {\it additive martingale},\nomenclature{$W_n$}{: the additive martingale}
\begin{equation}
\label{defW}
W_n:=\underset{|z|=n}{\sum}\ee^{-V(z)},
\end{equation}
and the probability measure $\Q$ such that for any $n\geq 0$,\nomenclature{$\Q$}{: law of the branching random walk with spine}
\begin{equation}
\Q|_{\mathcal{F}_n}:=W_n\bullet\P|_{\mathcal{F}_n},
\end{equation}
where $\mathcal{F}_n$ denotes the sigma-algebra generated by the positions $(V(z),\,|z|\leq n)$ up to time $n$. Let $\hat{L}$ be a point process whose has Radon-Nikodym derivative $\int\ee^{-x}L(dx)$ with respect to the law of $L$. In \cite{Lyo97}, Russell Lyons gave the following description of the branching random walk under $\Q$: 

-start with one particle $w_0$ at the origin. Generate offspring and displacements according to a copy $\hat{L}_1$ of $\hat{L}$,

-choose $w_1$ among children of $w_0$ with probability proportional to $\ee^{-(V(x))}$ when its displacement is $V(x)$,

-the children other than $w_1$ give rise to ordinary independent branching random walks,

-$w_1$ gives birth to particles distributed according to $\hat{L}$,

-again, choose one of the children of $w_1$ at random, call it $w_2$, with the others giving rise to ordinary independent branching random walks, and so on.

We still call $\mathbb{T}$ the genealogical tree of the process, so that $(w_n)_{n\in \N}$ is a ray of $\mathbb{T}$, which we will call the {\it spine}. This change of probability was also used in \cite{HuS09}. We refer to \cite{LPP95} for the case of the Galton-Watson tree, to \cite{CRo88} for the analogue for the branching Brownian motion, and to \cite{BKy04} for the spine decomposition in various types of branching.\nomenclature{$(w_n)_{n\in \N} $}{: spine of the branching random walk under $\Q$}
\begin{Proposition}
\label{lyons}
Under (\ref{boundarycase}),

(i) for any $|z|=n$, we have
\begin{equation}
\label{no(iii)}
\Q\left\{w_n=z\Big| \mathcal{F}_n\right\}=\frac{\ee^{-V(z)}}{W_n};
\end{equation}

(ii) the spine process $(V(w_n),n\geq 0)$ has distribution of the centered random walk $(S_n,n\geq 0)$ under $\Q$ satisfying (\ref{manytoone}).
\end{Proposition}

Before closing this subsection, we collect some useful facts about the centered random walks with finite variance $(S_n)_{n\in \N}$. We have taken these statements from Section 2 of \cite{Aid11}:
\begin{Lemma}
(i) There exists a constant $\alpha_1>0$ such that for any $x\geq 0$ and $n\geq 1$,
\begin{equation}
\P_x\left(\underset{j\leq n}{\min}\,S_j\geq 0\right)\leq \alpha_1(1+x)n^{-\frac{1}{2}}.
\label{2.5}
\end{equation}
\noindent{(ii) There} exists a constant $\alpha_2>0$ such that for any $b\geq a,x\geq 0$ and $ n\geq 1$ , 
\begin{equation}
\P_x\left(S_n\in[a,b],\underset{j\leq n}{\min}\,S_j\geq 0\right)\leq \alpha_2(1+x)(1+b-a)(1+b)n^{-\frac{3}{2}}.
\label{2.6}
\end{equation}
\noindent{(iii) Let} $0<\Lambda<1$. There exists a constant $\alpha_3=\alpha_3(\Lambda)>0$ such that for any $b\geq a ,x\geq 0, y\in \r$
\begin{eqnarray}
\label{2.7}
&&\P_x\left(S_n\in [y+a,y+b],\underset{j\leq n}{\min}\, S_j\geq 0, \underset{\Lambda n\leq j\leq n}{\min}\,S_j\geq y\right)\qquad\qquad\qquad\qquad
\\
&&\qquad\qquad\qquad\qquad\qquad\nonumber\leq \alpha_3(1+x)(1+b-a)(1+b)n^{-\frac{3}{2}}.
\end{eqnarray}
\end{Lemma}
\noindent See \cite{Koz76} for  (\ref{2.5}). The estimates   (\ref{2.6}) and (\ref{2.7}) are for example Lemmas A.1 and A.3 in \cite{AShi10}. (In our case $(S_n)$ is the centered random walk under $\P$, with finite variance $\E[S_1^2]=\sigma$ which appears in the many-to-one formula) 

We introduce its renewal function $R(x)$ which is zero if $x<0$, 1 if $x=0$, and \nomenclature{$R(x)$}{: renewal function associated to $(S)_{n\in \N}$}%
\begin{equation}
R(x):=\underset{k\geq 0}{\sum}\P\left(S_k\geq -x,\, S_k<\underset{0\leq j\leq k-1}{\min}S_j\right),\qquad \text{for }x>0.
 \label{2.9}
\end{equation}
It is known that there exists $c_0>0$ such that \nomenclature{$c_0$}{: defined in (\ref{2.10})}
\begin{equation}
 \underset{x\to\infty}{\lim}\frac{R(x)}{x}=c_0.
\label{2.10}
\end{equation}

Similarly, we define $R_-(x)$ as the renewal function associated to $-S$. Finally according to Theorem 1a, Section XII.7 p.415 of \cite{Fel71}, there exists $C_-,C_+>0$ such that 
\begin{equation}
\label{def2C}
\P\left( \underset{1\leq i\leq n}{\min }S_i\geq 0 \right) \sim \frac{C_+}{\sqrt{n}},\quad \P\left( \underset{1\leq i\leq n}{\max }S_i\leq 0 \right) \sim \frac{C_+}{\sqrt{n}},\quad \text{as } n\to \infty.
\end{equation}

\subsection{Notations}
We will use the notations of Aïdékon in \cite{Aid11}. As the typical order of $M_n$ is $\frac{3}{2}\log n$, it will be convenient to use the following notation, for $x\geq 0$:\nomenclature{$a_n(x) $, $I_n(x) $}{: $\frac{3}{2}\log n$ and $[a_n(x)-1,a_n(x))$}
\begin{eqnarray*}
\text{                    }&&a_n(x):=\frac{3}{2}\log n-x,
\\
\text{                    }&&I_n(x):=[a_n(x)-1,a_n(x)).
\end{eqnarray*}
Let us introduce for any $x_1,x_2,x_3>0$ the set $\mathcal{Z}^{x_1,x_2,x_3}$ defined by
$$ z\in \mathcal{Z}_n^{x_1,x_2,x_3} \Longleftrightarrow |z|=n,\,\underset{k\leq n}{\min}V(z_k)\geq -x_1,\, \underset{k\in [\frac{n}{2},n]}{\min} V(z_k)\geq a_n-x_2,\,  V(z)\leq a_n-x_3.$$ 

\subsection{On the tail of distribution of $W_{n,\beta}$.}
In this section we will study the tail of distribution of $W_{n,\beta}$. As we will see, when $A$ is large enough, up to a negligible amount, $W_{n,\beta} $ is equal to $ \sum_{|z|=n}\ee^{-\beta V(z)}\1_{\{V(z)\leq M_n +A \}}$. So we have to study the extremal particles, i.e the particles $ z\in \mathbb{T}_n$ such that $V(z)$ such that $M_n\simeq V(z)$. 

Let us recall two known results. For any $x\geq 0$,
\begin{eqnarray}
\nonumber \P\left(\exists u\in \mathbb{T}: V(u)\leq -x\right) &\leq& \sum_{n\geq 0}\E\Big( \sum_{|u|=n}\1_{\{V(u)\leq -x,\, V(u_k)\geq -x,\, \forall k<n \}} \Big)
\\
\nonumber &=& \sum_{n\geq 0}\E\big( \ee^{S_n},\, S_n\leq -x,\, S_k>-x,\ ,\forall k<n\big)
\\
\label{minoration} &\leq &\ee^{-x}\sum_{n\geq 0}\P\big(  S_n\leq -x,\, S_k>-x,\ ,\forall k<n\big) \leq \ee^{-x}.
\end{eqnarray}
where we have used (\ref{manytoone}) in the equality. From Corollary 3.4 \cite{Aid11} and (\ref{minoration}) it is plain to deduce:
\begin{Proposition}[\cite{Aid11}]
\label{tenstionaid}
Under (\ref{boundarycase}) and (\ref{condint}), there exists $c_{2}>0$ such that for any $x\geq 0$ and any integer $n\geq 1$,
\begin{equation}
\label{aidtension}
\P(M_n\leq a_n-x)\leq c_{2}(1+x)\ee^{-x}.
\end{equation}
\end{Proposition}
In this section we will principally show three results:
\begin{Proposition}
\label{A,L}
Under (\ref{boundarycase}) and (\ref{condint}), for any $\epsilon,K>0$ there exists $L_0$ large enough such that for any $L>L_0$ $n\in \N$, $x\in \r$, $y>0$ and $\delta\in [-\infty,K]$,

\begin{equation}
\label{eqpropA,L}
\P\left(n^{\frac{3}{2}\beta}\underset{|z|=n}{\sum}\ee^{-\beta V(z)}\1_{\{ z\notin \mathcal{Z}_n^{y,x+L,x-L}\}}\geq \ee^{\beta (x-\delta)}\right)\leq \epsilon (1+y)\ee^{-x}+\ee^{-y}.
\end{equation}
\end{Proposition}

\begin{Proposition}
\label{lemdomination2}
Under (\ref{boundarycase}) and (\ref{condint}), there exist $c_{3},\, c_4>0$ such that for any $n\in \N^*$, $j,\, x\geq 1$,
\begin{equation}
\label{eqlemdomination2tristris}
\P\left(\widetilde{W}_{n,\beta}\geq \ee^{\beta x},\, M_n \in I_n(x-j)\right)\leq c_{3} (1+x)\ee^{-x}\ee^{-c_4 j}.
\end{equation}
In particular we see that 
\begin{equation}
\label{eqlemdomination2}
\P(\widetilde{W}_{n,\beta}\geq \ee^{\beta x})\leq c_5 x\ee^{-x},\quad \forall n>1,\, x\geq 1.
\end{equation}
\end{Proposition}
And finally,
 \begin{corollary}
\label{partieneg}
(i) For any $ \epsilon>0$, there exists $L,\,  A,\, N>0$ large enough we have for any  $n\geq N$ and $x\geq A$, 
\begin{equation}
\label{eqpartieneg} \E\left(1- \exp\{-  \sum_{|z|=n}\ee^{-\beta [V(z)-a_n+x]} \1_{\{ V(z)\geq a_n-x+L \}} \} \right)\leq \epsilon x\ee^{-x} .
\end{equation}

\label{uplowbound}
(ii) There exist $c_6,\, c_7>0$ such that for any $ A,\, N>0$ large enough we have for any  $n\geq N$ and $x\in [A, \frac{3}{2}\log n-A]$, 
\begin{equation}
\label{equplowbound} c_6(1+x)\ee^{-x}\leq   \E\left(1- \exp\{-  \sum_{|z|=n}\ee^{-\beta [V(z)-a_n+x]}   \} \right)\leq c_7 (1+x)\ee^{-x}  .
\end{equation}
\end{corollary}
\paragraph{Remark:} Only (\ref{eqpartieneg}) and (\ref{equplowbound}) will be used in the next section.\\

\noindent{\it Proof of Corollary \ref{partieneg}.}
The lower bound of (\ref{equplowbound}) stems directly from
\begin{eqnarray*}
\E\left(1- \exp\{-  \sum_{|z|=n}\ee^{-\beta [V(z)-a_n+x]}   \} \right)&\geq& (1-\ee^{-1})\P\left(M_n\leq a_n-x\right)
\\
&\geq & c_6(1-\ee^{-1}) (1+x)\ee^{-x},
\end{eqnarray*}
where in the last line, we have used Proposition 4.1 in \cite{Aid11}. For the upper bound of (\ref{equplowbound}), notice that for any $w >0$, $1-\ee^{-w}= \int_0^{+\infty} \ee^{-x}\1_{\{w \geq x\}}dx$, thus according to Proposition \ref{lemdomination}, we can write
\begin{eqnarray*}
 \E\left(1- \exp\{-  \sum_{|z|=n}\ee^{-\beta [V(z)-a_n+x]}   \} \right)&=& \int_{0}^{+\infty}\ee^{-u} \P\left(\widetilde{W}_{n,\beta}\geq \ee^{\beta x}u \right)du
 \\
 &\leq & \ee^{-\beta (x-1)}+ c_1 \int_{\ee^{-\beta (x-1)}}^{+\infty} \frac{\ee^{-u}}{u^{\frac{1}{\beta}}} (x+ \frac{1}{\beta}\log u)\ee^{-x}du
 \\
 &\leq &  c_7 (1+x) \ee^{-x},
\end{eqnarray*}
which proves (\ref{equplowbound}). By using (\ref{eqpropA,L}) instead of (\ref{eqlemdomination2}), we derive (\ref{eqpartieneg}) identically; indeed it suffices to observe that $V(u)> a_n-x+L $ implies $ u\notin \mathcal{Z}^{x,x+L,x-L}$.
\hfill$\Box$
\\

Now we will prove Proposition \ref{A,L} and \ref{lemdomination}. First recall the important following Lemma:
\begin{Lemma}[Aïdékon \cite{Aid11}]
\label{afro}
 Under (\ref{boundarycase}) and (\ref{condint}), there exist constants $c_{8},c_{9}>0$ such that for any $n\geq 1,\,L\geq 0$ and $y\geq 0,\,x\in \r$,
\begin{equation}
\label{eqt4bis}
\P\left(\exists u\in \mathbb{T}_n,\,\underset{k\leq n}{\min}{V}(u_k)\geq -y  ,\,\underset{\frac{n}{2}<k\leq n}{\min}V(u_k)\leq a_n(x+L),\,V(u)\leq a_n(x)   \right) \leq c_{8}(1+y)\ee^{-c_{9}L}\ee^{-x}.
\end{equation}
\end{Lemma}
(We have given a statement slightly stronger than in \cite{Aid11} but reader can see easily that actually this is an equivalent statement).

The proofs require two Lemmas. 

\begin{Lemma}
Under (\ref{boundarycase}) and (\ref{condint}), there exist $c_{10},c_{11}>0$ such that $\forall x\in \r,\,y,\,L\geq 0$, $ n\geq 1$,
\begin{equation}
\label{eqt3bis}
\P_y\left(n^{\frac{3}{2}}\underset{|z|=n}{\sum}\ee^{-\beta V(z)}\1_{\{\underset{k\leq n}{\min}V(z_k)\geq 0,\underset{\frac{n}{2}<k\leq n}{\min}V(z_k)\leq a_n(x+L)\}}\geq  \ee^{\beta x} \right)\leq c_{10}(1+y)\ee^{-c_{11} L}\ee^{-y} \ee^{-x}.
\end{equation}
\label{t3bis}
\end{Lemma}
\noindent{\it Proof of Lemma \ref{t3bis}.}
For any $a,b>0$ let $\P_{(\ref{eqt4bis})}(a,b)$ and $\P_{(\ref{eqt3bis})}(a,b)$ be the probability of respectively (\ref{eqt4bis}) and (\ref{eqt3bis}) when $y=a$ and $x=b$.  Observe that
\begin{eqnarray}
\nonumber \P_{(\ref{eqt3bis})}(y,x) &\leq& \P_y\left(\underset{|z|=n}{\sum}\ee^{-\beta V(z)}\1_{\{\underset{k\leq n}{\min}V(z_k)\geq 0,\underset{\frac{n}{2}<k\leq n}{\min}V(z_k)\leq a_n(x+L),V(z)\geq a_n(x)\}})\geq \frac{\ee^{\beta x}}{n^{\frac{3}{2}\beta}}\right)
\\
&&\nonumber  \qquad\qquad\qquad\qquad\qquad\qquad  \qquad\qquad\qquad \qquad\qquad\qquad  +\P_{(\ref{eqt4bis})}(y,y+x)
\\
\label{quatredollars}&  \leq& \P_y(...)+c_{8}(1+y)\ee^{-c_{9}L}\ee^{-x-y}.
\end{eqnarray}
We have to bound $\P_y(...)$.

We need some notations: for $|z|=n$, $j\geq 0$, $\frac{n}{2}<k\leq n$ and $L'\geq L$ we define the event 
\begin{equation}
\label{44.166}
E_{j,k,L'}(z):=\{\underset{l\leq n}{\min}V(z_l)\geq 0,\,V(z_k)=\underset{\frac{n}{2}<l\leq n}{\min}V(z_l)\in I_n(x+L'),\, V(z_n)\in I_n(x)+j\}.
\end{equation}
For any integer $a\in [0,\frac{n}{2}]$ let\nomenclature{$E_{j,k,L'}(z),\,F_{L'}^{[\frac{n}{2},n-a)}(z),\,F_{L'}^{[n-a,n]}(z),\,E_{j,k,L'}(S),\,F_{L'}^{[\frac{n}{2},n-a)}(S),\,F_{L'}^{[n-a,n]}(S) $}{: defined in (\ref{44.166}) to (\ref{44.199}) }
\begin{equation}
F_{L'}^{[\frac{n}{2},n-a]}(z):=\underset{j\geq 0,\,k\in[\frac{n}{2},n-a]}{\bigcup}E_{j,k,L'}(z),\qquad F^{(n-a,n]}_{L'}(z):=\underset{j\geq 0,\,k\in(n-a,n]}{\bigcup}E_{j,k,L'}(z).
\end{equation}
Remark that for any sequence of integer $(a_{L+p})_{p\geq 0}$, 
\begin{eqnarray}
&&\nonumber \underset{|z|=n}{\sum}\ee^{-\beta V(z)}\1_{\{\underset{k\leq n}{\min}V(z_k)\geq 0,\underset{\frac{n}{2}<k\leq n}{\min}V(z_k)\leq a_n(x+L),V(z)\geq a_n(x)\}}
\\
&&\label{deuxdollar}\qquad\qquad\qquad \qquad \leq \underset{p\geq 0}{\sum}\underset{|z|=n}{\sum}\ee^{-\beta V(z)}
(\1_{\{F^{[\frac{n}{2},n-a_{L+p}]}_{L+p}(z)\}} + \1_{\{F^{(n-a_{L+p,n}]}_{L+p}(z)\}}) .
\end{eqnarray}

Similarly we introduce, for the centered random walk $(S_n)_{n\geq 0}$,
\begin{eqnarray}
E_{j,k,L'}:=\{\underset{l\leq n}{\min}\,S_l\geq 0,\,S_k=\underset{\frac{n}{2}<l\leq n}{\min}S_l\in I_n(x+L'),\, S_n\in I_n(x)+j\},
\\
 \qquad F^{[\frac{n}{2},n-a]}_{L'}:=\underset{j\geq 0,\,k\in[\frac{n}{2},n-a]}{\bigcup}E_{j,k,L'},\qquad F^{(n-a,n]}_{L'}:=\underset{j\geq0,\,k\in(n-a,n]}{\bigcup}E_{j,k,L'}.\label{44.199}
\end{eqnarray}

Let estimate $\P_y(E_{j,k,L'})$ for $\frac{n}{2}< k\leq n-a$. By the Markov property at time k,
\begin{eqnarray*}
\P_y(E_{j,k,L'})\leq \P_y\left(\underset{l\leq k}{\min} \,S_l\geq 0,\, \underset{\frac{n}{2}<l\leq k}{\min} S_l\geq a_n(x+L'),\, S_k\in I_n(x+L')\right)\times 
\\
\P\left(S_{n-k}\in[L'-1+j,L'+1+j],\underset{l\leq n-k}{\min}S_l\geq 0\right).
\end{eqnarray*}
We know by (\ref{2.6}) that there exists a constant $c_{12}$ such that
$$\P\left(S_{n-k}\in[L'-1+j,L'+1+j],\,\underset{l\leq n-k}{\min}\,S_l\geq 0\right)\leq c_{12}(n-k+1)^{-\frac{3}{2}}(1+L'+j).$$
For the first term, we have to discuss on the value of $k$. Suppose that $\frac{3}{4}n\leq k\leq n$, then by (\ref{2.7}),
$$\P_y\left(\underset{l\leq k}{\min}\, S_l\geq 0,\, \underset{\frac{n}{2}<l\leq k}{\min} S_l\geq a_n(x+L'),\, S_k\in I_n(x+L') \right)\leq c_{13}\frac{1+y}{n^{\frac{3}{2}}}.$$
If $\frac{1}{2}n\leq k\leq \frac{3}{4}n$ we simply write
\begin{eqnarray*}
\P_y\left(\underset{l\leq k}{\min}\, S_l\geq 0,\, \underset{\frac{n}{2}<l\leq k}{\min} S_l\geq a_n(x+L'),\, S_k\in I_n(x+L')\right)&\leq& \P_y\left(S_k\in I_n(x+L'),\underset{l\leq k}{\min} S_l\geq 0\right)
\\
&\leq& c_{14}(1+y)n^{-\frac{3}{2}}\log n .
\end{eqnarray*}
To summarize we have obtained
\begin{equation}
\P_y(E_{j,k,L'})\leq\left\{
         \begin{array}{ll}
           c_{15}\frac{(1+y)\log n}{n^{\frac{3}{2}}(n-k+1)^{\frac{3}{2}}}(1+L'+j) & \mathrm{if}\quad \frac{n}{2}<k\leq \frac{3}{4}n, \\
            c_{15}\frac{1+y}{n^{\frac{3}{2}}(n-k+1)^{\frac{3}{2}}}(1+L'+j)& \text{if} \qquad \frac{3}{4}n<k\leq n-a. \\
         \end{array}
       \right.
\end{equation}
Now we can tackle the proof and study $\P_y(...)$. By Proposition \ref{lyons},
\begin{eqnarray*}
\frac{n^{\frac{3}{2}\beta}}{\ee^{\beta x}} \E_y\left(\underset{|z|=n}{\sum}\ee^{-\beta V(z)}
\1_{\{F^{[\frac{n}{2},n-a]}_{L'}(z)\}}\right)&=&\frac{n^{\frac{3}{2}\beta}}{\ee^{\beta x}}\ee^{-y}\E_y\left[\ee^{(1-\beta) S_n}\1_{\{F^{[\frac{n}{2},n-a]}_{L'}\}}\right]
\\
&\leq& \frac{n^{\frac{3}{2}\beta}}{\ee^{\beta x}}\ee^{-y}\overset{n-a}{\underset{k=n/2}{\sum}} \underset{j\geq 0}{\sum} \ee^{(1-\beta)(a_n(x)+j)}\P_y(E_{j,k,L'})
\\
&\leq& c_{16} (1+y)(1+L')\ee^{-x-y}a^{-\frac{1}{2}}.
\label{t0}
\end{eqnarray*}
We also get
\begin{eqnarray*}
&&\P_y\left(\underset{|z|=n}{\sum}\1_{\{F^{(n-a,n]}_{L'}(z)\}}\geq 1\right)\leq\underset{k\in[n-a,n]}{\sum}\P_y\left(\underset{|z|=n}{\sum}\underset{j\geq 0}{\sum}\1_{\{E_{j,k,L'}(z)\}}\geq 1\right)
\\
&&\leq \underset{k\in(n-a,n]}{\sum} \P_y\left(\exists |z|=k: \underset{l\leq k}{\min}\, V(z_l)\geq 0,\, \underset{\frac{n}{2}<l\leq k}{\min} V(z_l)\geq a_n(x+L'),\, V(z_k)\in I_n(x+L')\right).
\end{eqnarray*}
Using Proposition \ref{lyons} then (\ref{2.7}), we deduce that 
\begin{eqnarray*}
\P_y\left(\underset{|z|=n}{\sum}\1_{\{F^{(n-a,n]}_{L'}(z)\}}\geq 1\right)\leq  \underset{k\in[n-a,n]}{\sum} c_{17}(1+y)\ee^{-x-y-L'}= c_{17}(1+a)(1+y)\ee^{-x-y-L'}.
\end{eqnarray*}

These estimations lead us to choose correctly the integer $a$. Let $\alpha\in (0,1)$, for any $p\geq 0$ we define $a_{L+p}:=\lfloor\ee^{\alpha(L+p)}\rfloor$ . Recalling (\ref{deuxdollar}), to conclude we just assemble our previous inequalities and observe that
\begin{eqnarray}
\nonumber \P_y{(...)}&\leq &\P_y\left(\underset{p\geq 0}{\sum}\underset{|z|=n}{\sum}\ee^{-\beta V(z)}
\1_{\{F^{[\frac{n}{2},n-a_{L+p}]}_{L+p}(z)\}}\geq \frac{\ee^{\beta x}}{2n^{\frac{3}{2}\beta}}\right)+\P_y\left(\underset{p\geq 0}{\sum}\underset{|z|=n}{\sum}\1_{\{F^{(n-a_{L+p},n]}_{L+p}(z)\}}\geq 1\right)
\\
\nonumber&\leq&\underset{p\geq 0}{\sum}\left(\frac{n^{\frac{3}{2}\beta}}{\ee^{\beta x}} \E_y\left(\underset{|z|=n}{\sum}\ee^{-\beta V(z)}
\1_{\{F^{[\frac{n}{2},n-a_{L+p}]}_{L+p}(z)\}}\right)+ \P_y\left(\underset{|z|=n}{\sum}\1_{\{F^{(n-a_{L+p},n]}_{L+p}(z)\}}\geq 1\right)\right)
\\
\nonumber&\leq&\underset{p\geq 0}{\sum}\left(c_{16} (1+y)(1+L+p)\ee^{-x-y}a_{L+p}^{-\frac{1}{2}}+c_{17}(1+a_{L+p})(1+y)\ee^{-x-y-(L+p)}\right)
\\
\label{troisdollars} &\leq& c_{10}(1+y)\ee^{-c_{11} L}\ee^{-y-x}.
\end{eqnarray}
Combining (\ref{quatredollars}) with (\ref{troisdollars}) we obtain Lemma \ref{t3bis}. \hfill $\Box$

\begin{Lemma} 
Under (\ref{boundarycase}) and (\ref{condint}), there exist $c_{18}\geq 0$, $c_{19}>0 $ ($c_{19}={\beta -1}$) such that for any $n\geq 0$, $x\in \r$, $A,\, y\geq 0$ and $L\geq 0$,
\begin{eqnarray*}
\P_y\Big(\underset{|z|=n}{\sum}\ee^{-\beta V(z)}\1_{\{\underset{k\leq n}{\min}V(z_k)\geq 0,\underset{\frac{n}{2}<k\leq n}{\min}V(z_k)\geq a_n(x+L),\,V(z)\geq a_n(x)+A\}}\geq \frac{\ee^{\beta x}}{n^{\frac{3}{2}\beta}}\Big) \qquad\qquad 
\\
\leq c_{18}(1+y)\ee^{-x-y}(L+A+1)\ee^{-c_{19} A}.
\end{eqnarray*}
\label{lemme3}
\end{Lemma}
\noindent{\it Proof of Lemma \ref{lemme3}.} The Markov inequality $\P(X\geq 1)\leq \E(X)$, for $X$ positive gives:
\begin{eqnarray*}
&&\P_y\Big(\underset{|z|=n}{\sum}\ee^{-\beta V(z)}\1_{\{\underset{k\leq n}{\min}V(z_k)\geq 0,\,\underset{\frac{n}{2}<k\leq n}{\min}V(z_k)\geq a_n(x+L),\,V(z)\geq a_n(x)+A\}}\geq \frac{\ee^{\beta x}}{n^{\frac{3}{2}\beta}}\Big)\leq
\\
&& \frac{n^{\frac{3}{2}\beta}}{\ee^{\beta x}}\E_y\Big(\underset{|z|=n}{\sum}\ee^{-\beta V(z)}\1_{\{\underset{k\leq n}{\min}V(z_k)\geq 0,\,\underset{\frac{n}{2}<k\leq n}{\min}V(z_k)\geq a_n(x+L),\,V(z)\geq a_n(x)+A\}}\Big).
\end{eqnarray*}
By Proposition \ref{lyons} this is equal to
\begin{eqnarray*}
&=& \frac{n^{\frac{3}{2}\beta }}{\ee^{\beta x}}\ee^{-y}\underset{k\in \N}{\sum}\E_y\left(\ee^{(1-\beta) S_n}\1_{\{\underset{k\leq n}{\min}\,S_k\geq 0,\,\underset{\frac{n}{2}<k\leq n}{\min}S_k\geq a_n(x+L),\, S_n\in I_n(x-A-k)\}}\right)
\\
&\leq& \ee^{-x-y-A(\beta-1)}\underset{k\in \N}{\sum}\ee^{(1-\beta)k}n^{\frac{3}{2}}\P_y\left(\underset{k\leq n}{\min}\,S_k\geq 0,\,\underset{\frac{n}{2}<k\leq n}{\min}S_k\geq a_n(x+L),\, S_n\in I_n(x-A-k)\right)
\\
&\leq & \ee^{-y-x}\ee^{-A(\beta-1)} \underset{k\in \N}{\sum}\ee^{(1-\beta)k}c_{20}(1+y)(L+A+k+1)
\\
&\leq & c_{18}(1+y)\ee^{-x-y}(L+A+1)\ee^{-c_{19} A},
\end{eqnarray*}
where in the second inequality we have used used (\ref{2.7}).
\hfill $\Box$
\\

\noindent{\it Proof of Proposition \ref{A,L}.} By noticing that $ \mathcal{Z}_n^{y,x+L,x-L} =\mathcal{Z}_n^{y,x-\delta+L+\delta,x-\delta-(L-\delta)}\subset \mathcal{Z}_n^{y,x-\delta+ L+K,x-\delta-(L+K)}$, it is sufficient to prove the Proposition \ref{A,L} for $K=0$. Let $z\in \mathbb{T}_n$, recall that
\begin{equation}
z\in \mathcal{Z}^{y,x+L,x-L}_n\Longleftrightarrow \underset{k\leq n}{\min}{V}(z_k)\geq -y,\, \underset{k\in [n/2,n]}{\min}{V}(z_k)\geq a_n-(x+L),\, V(z)\leq a_n-(x-L).
\end{equation}
As $\P\left(\exists z\in \mathbb{T},\, V(z)\leq -y\right)\leq \ee^{-y}$ we deduce that the probability (\ref{eqpropA,L}) is smaller than
\begin{equation}
\P\Big(n^{\frac{3}{2}\beta}\sum_{|z|=n}\ee^{-\beta V(z)} \1_{\{ \underset{k\leq n}{\min}{V}(z_k)\geq -y,\, z\notin \mathcal{Z}_n^{y,x+L,x-L}\}} \geq \ee^{\beta x} \Big)+\ee^{-y}.
\end{equation}
 Then we observe that on $\{ \underset{k\leq n}{\min}V(z_k)\geq -y\}$,
{\small \begin{eqnarray*}
&&\{z\notin \mathcal{Z}^{y,x+L,x-L}\}\subset
 \{\underset{\frac{n}{2}<k\leq n}{\min}V(z_k)\leq a_n(x+L)\}\cup \{\underset{\frac{n}{2}<k\leq n}{\min}V(z_k)\geq a_n(x+L),\,V(z)\geq a_n(x)+L\}.
\end{eqnarray*}}
According to Lemmas \ref{t3bis} and \ref{lemme3}, by choosing $L_0$ large enough such that $\forall L\geq L_0$, $c_{10}\ee^{-c_{11}L}+c_{18}(2L+1)\ee^{-c_{19} L}\leq \epsilon$, we obtain for any $n\in \N$, $x\in \r$ and $y>0$, 
\begin{equation}
\P\left(n^{\frac{3}{2}\beta}\underset{|z|=n}{\sum}\ee^{-\beta V(z)}\1_{\{ z\notin \mathcal{Z}_n^{y,x+L,x-L}\}}\geq \ee^{\beta x}\right)\leq \epsilon (1+y)\ee^{-x}+\ee^{-y},
\end{equation}
which proves (\ref{eqpropA,L}).
\hfill$\Box$
\\
\\

\noindent{\it Proof of Proposition \ref{lemdomination2}.} If $M_n\geq a_n(x)+j$, then clearly for any $z\in \mathbb{T}$ such that $|z|=n$ we have $z\notin \mathcal{Z}^{x,x+j,x-j}$. Therefore,
\begin{eqnarray*}
\P\left(\widetilde{W}_{n,\beta}\geq \ee^{\beta x}, M_n\in I_n(x-j+1)\right) \leq \P\left(n^{\frac{3}{2}\beta}\underset{|z|=n}{\sum}\ee^{-\beta V(z)}\1_{\{ z\notin \mathcal{Z}^{x+j,x+j,x-j}\}}\geq  {\ee^{\beta x}} \right)\qquad\qquad\quad 
\\
\leq  c_{10}(1+x+j)\ee^{-x}\ee^{-c_{11} j}+  c_{18}(1+x+j)\ee^{-x}j\ee^{-c_{19} j}+\ee^{-(x+j)},
\end{eqnarray*}
by Lemma \ref{t3bis}, Lemma \ref{lemme3} and (\ref{minoration}). Set $c_{11}=\frac{1}{2}\min(c_{14},c_{22})$ to obtain (\ref{eqlemdomination2tristris}). The estimate (\ref{eqlemdomination2}) follows easily from (\ref{eqlemdomination2tristris}) and Proposition \ref{tenstionaid} applied to the equality
$$ \P\left(\widetilde{W}_{n,\beta}\geq \ee^{\beta x}\right)= \P\left(\widetilde{W}_{n,\beta}\geq \ee^{\beta x},\,M_n\leq a_n(x)\right)+\P\left(\widetilde{W}_{n,\beta}\geq \ee^{\beta x},\, M_n> a_n(x)\right). $$
This proves Proposition \ref{lemdomination2}.
\hfill $\Box$

\section{Proof of Proposition \ref{Propconclusion2}}
We fix $d\geq 1$, ${\bm \beta}\in (1,\infty)^d$, ${\bm \theta }\in (\r_{+}^*)^d$ and $\epsilon>0$. Our aim is to study for $n$ then $x$ large,
\begin{equation}
\E\Big(  1-\exp\{ -\sum_{i=1}^{d}   \sum_{|z|=n}\theta_i\ee^{-\beta_i[V(z)+x-a_n]} \}\Big).
\end{equation}
According to  (\ref{eqpartieneg}), for $L$ large enough we can restrain our study to the expectation of 
\begin{equation}
\Phi^{(L)}(x,n):= 1-\exp\{\sum_{i=1}^{d}   \sum_{|z|=n}\theta_i\ee^{-\beta_i[V(z)+x-a_n]}\1_{\{ V(z)\leq a_n-x+L \}} \} .
\end{equation}
Indeed as for any $u,v\geq 0$,  $1-\ee^{-u-v}\leq 1-\ee^{-u}+1-\ee^{-v}$, we deduce that
\begin{eqnarray*}
&&\left|\E\Big(  1-\exp\{ -\sum_{i=1}^{d}   \sum_{|z|=n}\theta_i\ee^{-\beta_i[V(z)+x-a_n]} \}\Big)-\E(\Phi^{(L)}(x,n))\right|
\\
&&\leq \E\Big( 1-\exp\{\sum_{i=1}^{d}   \sum_{|z|=n}\theta_i\ee^{-\beta_i[V(z)+x-a_n]}\1_{\{ V(z)\geq a_n-x+L \}} \} \Big).
\end{eqnarray*}

On the set $\{ M_n > a_n-x+L\}$, $\Phi^{(L)}(x,n)=0$, thus we have
\begin{eqnarray}
\nonumber  \E(\Phi^{(L)}(x,n))&=& \E\left(1-\exp\{-\sum_{i=1}^{d}   \sum_{|z|=n}\theta_i\ee^{-\beta_i[V(z)+x-a_n]}\1_{\{ V(z)\leq a_n-x+L \}} \};  M_n \leq a_n-x+L  \right)
\\
\label{5.5}&=&  \E\left(\sum_{|z|=n}  \frac{  \1_{\{ V(z)=M_n \leq a_n-x+L \}}  }{\sum_{|z|=n}\1_{\{ V(z)=M_n\}}} \Phi^{(L)}(x,n)  \right)  :=\E_{(\ref{5.5})}.
\end{eqnarray}
According to the Lemma \ref{afro} and (\ref{minoration}), we can enhance $L$ such that for any $n,x\geq 0$ large enough,
\begin{equation}
\P\left(M_n \leq a_n-x+L,\, \exists z\notin \mathcal{Z}^{x,x+ \frac{2}{c_9} L,x-L} \right)\leq \epsilon x\ee^{-x},
\end{equation}
whrere $c_9$ is the constant defined in (\ref{eqt4bis}). The random variable in (\ref{5.5}) is smaller than 1, then we deduce that for any $x,\, n\geq 1$ large enough, 
\begin{equation}
\label{xan} |\E_{(\ref{5.5})}-\E\left(\sum_{|z|=n}  \frac{  \1_{\{ V(z)=M_n\leq a_n-x+L,\, z\in \mathcal{Z}^{x,x+ \frac{2}{c_9}L,x-L} \}}  }{\sum_{|z|=n}\1_{\{ V(z)=M_n\}}} \Phi^{(L)}(x,n)  \right)|\leq  \epsilon x\ee^{-x}.
\end{equation}
Combining (\ref{xan}) to
\begin{eqnarray*}
 \sum_{|z|=n}  \frac{  \1_{\{ V(z)=M_n\leq a_n-x+L,\, z\in \mathcal{Z}^{x,x+ \frac{2}{c_9}L,x-L} \}}  }{\sum_{|z|=n}\1_{\{ V(z)=M_n\}}} \Phi^{(L)}(x,n) \qquad\qquad\qquad
 \\
  =  \sum_{|z|=n}  \frac{  \1_{\{ V(z)=M_n\leq a_n-x+\frac{2}{c_9} L,\, z\in \mathcal{Z}^{x,x+ \frac{2}{c_9}L,x- \frac{2}{c_9}L} \}}  }{\sum_{|z|=n}\1_{\{ V(z)=M_n\}}} \Phi^{(L)}(x,n),
\end{eqnarray*}
and 
\begin{eqnarray*}
\E\left(|\Phi^{(L)}(x,n)-  \Phi^{(\frac{2}{c_9} L)}(x,n)|\right) \leq  \E\Big( 1-\exp\{\sum_{i=1}^{d}   \sum_{|z|=n}\theta_i\ee^{-\beta_i[V(z)+x-a_n]}\1_{\{ V(z)\geq a_n-x+L \}} \} \Big) \leq \epsilon x\ee^{-x},
\end{eqnarray*}
we finally deduce the following statement: {\it  there exists $L'(=\frac{2}{c_9}L)$ large enough such that for any $n,x$ large enough
\begin{eqnarray}
\nonumber \Big|\E\Big(  1-\exp\{ -\sum_{i=1}^{d}   \sum_{|z|=n}\theta_i\ee^{-\beta_i[V(z)+x-a_n]} \}\Big)-  \E\Big(\sum_{|z|=n}  \frac{  \1_{\{ V(z)=M_n\leq a_n-x+  L',\, z\in \mathcal{Z}^{x,x+  L',x- L'} \}}  }{\sum_{|z|=n}\1_{\{ V(z)=M_n\}}} \Phi^{(L')}(x,n) \Big) \Big| 
\\
\label{prioqsopdij} \leq \epsilon x\ee^{-x}.
\end{eqnarray}}
From now we fix $L$ large such (\ref{prioqsopdij}) is true and study,
\begin{eqnarray}
\nonumber  &&\E\Big(\sum_{|z|=n}  \frac{  \1_{\{ V(z)=M_n\leq a_n-x+  L,\, z\in \mathcal{Z}^{x,x+ L,x- L} \}}  }{\sum_{|z|=n}\1_{\{ V(z)=M_n\}}} \Phi^{(L)}(x,n) \Big)
\\
\label{5.26} &&= \E\left( \ee^{V(w_n)}\1_{\{V(w_n)=M_n, w_n\in \mathcal{Z}^{x,x+L,x-L}\}}\frac{\Phi^{(L)}(x,n)}{\sum_{|z|=n}\1_{\{ V(z)=M_n\}} }\right).
\end{eqnarray}
Note that we have used the Proposition \ref{lyons} in the equality (\ref{5.26}). We denote by $\E_{(\ref{5.26})}$ the expectation in (\ref{5.26}).
\begin{definition}
\label{deuzedef}
For $b$ integer, we define the event $\xi_n$ by\nomenclature{$\xi_n(x,b,A) $}{: defined in (\ref{defxxi}). Good event on which the particles at generation $n$ which are located below $a_n(z)+A$ have a common ancestor with the spine at generation greater than $n-b$}
\begin{equation}
\xi_n:=\xi_n(x,b,L):=\{\forall k\leq n-b, \forall v\in \Omega(w_k), \underset{u\geq v,|u|=n}{\min}V(u)>a_n-x+L\},
\label{defxxi}
\end{equation}
where $\Omega(w_k)$ denotes the set of brothers of $w_k$. On the event $\xi_n\cap\{M_n\leq a_n-x+L\}$ we are sure that any particle located at the minimum separated from the spine after the time $n-b$. 
\end{definition}
\begin{definition}
\label{seuldef}
For $y\in \r,\,L,>0,\, {\bm \theta}\in (\r_{+}^*)^d$ and $b\in \N^*$ we define

(i) the function  
\begin{eqnarray}
\label{defdeF} F_{L,b}({\bm \theta},y):=\E_{\Q_y}\Big[\frac{\ee^{V(\omega_b)-L}\1_{\{V(\omega_b)=M_b\}}}{\underset{|u|=b}{\sum}\1_{\{V(u)=M_b\}}}\1_{\{V(\omega_b)\leq 2L,\, \underset{k\leq b}{\min}{V}(\omega_k)\geq 0 \}}\times \qquad\qquad\qquad
\\
\nonumber  \big[  1-\exp\{-\sum_{i=1}^{d} \sum_{|z|=b}\theta_i\ee^{-\beta_i[V(z)+y-L]}\1_{\{ V(z)+y\leq 2L \}}\big]  \Big];
\end{eqnarray}
(iii) the constant $\rho^{*}_{{\bm \beta},L,b}({\bm \theta}):=\frac{C_-C_+\sqrt{\pi}}{\sigma\sqrt{\pi}}\int_{y\geq 0}F_{L,b}({\bm \theta},y)R_-(y)dy$, where $C_-, C_+$ and $R_-(x)$ are defined in (\ref{def2C}).
\end{definition}
The proof of the following Lemma (which is an extension of Lemma 3.8 in \cite{Aid11}) is postponed in the Appendix C.
\begin{Lemma}
\label{t17}
Under (\ref{boundarycase}) and (\ref{condint}), $\forall \eta,L>0$ $\exists D(L,\eta)>0$  and  $B(L,\eta) \geq 1$ such that $\forall b\geq B,\, n\geq \ee^{5b},\,x\geq D$,
\begin{equation}
\Q\left((\xi_n(b,x,L))^c,\,  w_n\in \mathcal{Z}_n^{x,x+L,x-L}\right)\leq \eta n^{-\frac{3}{2}} (1+x).
\label{t16bis}
\end{equation}
\end{Lemma}
Lemma \ref{t17} justifies the definitions \ref{deuzedef} and \ref{seuldef}. Indeed observe that on $\xi_n(b,x,L) $, $ {\sum_{|z|=n}\1_{\{V(z)=M_n \}}} =  {\sum_{|z|=n,z>w_{n-b}}\1_{\{V(z)=M_n \}}} $  and $\Phi^{(L)}(x,n)$ is equal to
\begin{eqnarray*}
&& \Phi^{(L,b)}(x,n,w_{n-b})=    1-   \exp\{-\sum_{i=1}^{d}    \sum_{|z|=n,z>w_{n-b}}\theta_i\ee^{-\beta_i[V(z)-a_n+x]}\1_{\{ V(z)\leq a_n-x+L \}} \} .
\end{eqnarray*}
The functions $\Phi^{(L,b)}(\cdot,\cdot,\cdot)$ and $\Phi^{(L)}(\cdot,\cdot)$ are bounded by $1$, then applying Lemma \ref{t17}, there exists $b>0$ large enough (associated to $L,\, \frac{\eta}{\ee^L}$) such that 
\begin{eqnarray*}
&&|\E\Big( \ee^{V(w_n)}\1_{\{V(w_n)=M_n, w_n\in \mathcal{Z}^{x,x+L,x-L}\}}\big(\frac{\Phi^{(L)}(x,n)}{\sum_{|z|=n}\1_{\{ V(z)=M_n\}} } - \frac{ \Phi^{(L,b)}(x,n,w_{n-b})}{\sum_{|z|=n,z>w_{n-b}}\1_{\{V(z)=M_n \}}}\big) \Big)|.
\\
&&\leq \E\Big( \ee^{V(w_n)}  \big[\frac{ \Phi^{(L,b)}(x,n,w_{n-b})}{\sum_{|z|=n,z>w_{n-b}}\1_{\{V(z)=M_n \}}}+  \frac{ \Phi^{(L)}(x,n)}{\sum_{|z|=n}\1_{\{V(z)=M_n \}}}\big];V(w_n)=M_n,  w_n\in \mathcal{Z}^{x,x-L,x+L},\, 
\\
 &&\qquad\qquad\qquad\qquad  \qquad\qquad\qquad   \qquad \qquad   \qquad\qquad\qquad   \qquad \qquad    \qquad   \qquad \qquad  (\xi_n(b,x,L))^c\Big)
\\
&&\leq  \ee^{-x+L}n^{\frac{3}{2}} \Q\left((\xi_n(b,x,L))^c,\,  w_n\in \mathcal{Z}_n^{x,x+L,x-L}\right)
\\
&&\leq \eta \ee^{-x} (1+x).
\end{eqnarray*}

Moreover, using the Definition \ref{seuldef} and the branching property at time $n-b$ we have
\begin{eqnarray*}
\nonumber &&\E\Big( \ee^{V(w_n)}\1_{\{V(w_n)=M_n, w_n\in \mathcal{Z}^{x,x+L,x-L}\}} \frac{ \Phi^{(L,b)}(x,n,w_{n-b})}{\sum_{|z|=n,z>w_{n-b}}\1_{\{V(z)=M_n \}}} \Big)
\\
\nonumber &&= {n^{\frac{3}{2}}}\ee^{-x}\E_\Q\left(\1_{\{\underset{k\leq n-b}{\min}{V}(w_{k})\geq -x,\, \underset{k\in [\frac{n}{2},n-b]}{\min}{V}(w_{k})\geq a_n-x-L\}} F_{L,b}({\bm \theta}, V(w_{n-b})-a_n+x+L)\right)
\\
&&= {n^{\frac{3}{2}}}\ee^{-x}\E_{\Q_{x}}\left(\1_{\{\underset{k\leq n-b}{\min}{V}(w_k)\geq 0,\, \underset{k\in [\frac{n}{2},n-b]}{\min}{V}(w_k)\geq a_n-L\}} F_{L,b}({\bm \theta}, V(w_{n-b})-a_n+L)\right).
\end{eqnarray*}
By combining this equality with (\ref{prioqsopdij}) and (\ref{5.26}), we can affirm that there exists $L,\, B,\, D\geq 0$ such that for any $b\geq B$, $n\geq \ee^{5b}$,  $x\geq D$, 
\begin{eqnarray}
\nonumber &&\left| \E\Big(  1-\exp\{-\sum_{i=1}^{d}   \sum_{|z|=n}\theta_i\ee^{-\beta_i[V(z)+x-a_n]} \}\Big)- \right.
\\
\label{petitpoint}&&\left. {n^{\frac{3}{2}}}\ee^{-x}\E_{\Q_{x}}\left(\1_{\{\underset{k\leq n-b}{\min}{V}(w_k)\geq 0,\, \underset{k\in [\frac{n}{2},n-b]}{\min}{V}(w_k)\geq a_n-L\}} F_{L,b}({\bm \theta}, V(w_{n-b})-a_n+L)\right)\right|\leq \epsilon x\ee^{-x}.
\end{eqnarray}

Keeping in mind this last display, we shall now state and prove yet two lemmas which will be used in the proof of the Proposition \ref{Propconclusion2} :

\begin{Lemma}
\label{Riemann}
For any ${\bm \theta}\in (\r_{+}^*)^d$, the function $y\mapsto F_{L,b}({\bm \theta},y)$ is Riemann integrable and there exists a non-increasing function $\bar{F}:\r_+\to\r$ such that $|F_{L,b}({\bm \theta},y)|\leq \bar{F}(y)$ for any $y\geq 0$ and $\int_{y\geq 0}y\bar{F}(y)dy<\infty$.
\end{Lemma}
\noindent{\it Proof of Lemma \ref{Riemann}.} 
We recall that by Proposition \ref{lyons} the spine has the law of $(S_n)_{n\geq 0}$. We see that $\frac{\1_{\{V(\omega_b)=M_b\}}}{\underset{|u|=b}{\sum}\1_{\{V(u)=M_b\}}}$ is smaller than $1$, and $\ee^{V(\omega_b)-L}\leq \ee^{L}$. Hence, $|F_{L,b}({\bm \theta},y)|\leq \ee^{L} \P(S_b\leq L-y)=:\overset{-}{F}(y)$ which is non-increasing in $y$, and $\int_{y\geq 0}\overset{-}{F}(y)y dy= \ee^{L}\frac{1}{2}\E[(L-S_b)^21_{\{S_b\leq L\}}]<\infty$. Moreover, using the identity $|\1_E-a\1_F|\leq 1-a+|\1_E-\1_F|$ for $a\in(0,1)$, it yields that for $y_2\geq 0$,$\epsilon>0$ and any $y_1\in[y_2,y_2+\epsilon]$,
\begin{eqnarray*}
\left|F_{L,b}({\bm \theta},y_1)-F_{L,b}({\bm \theta},y_2)\right|\leq \E_\Q\left[|\1_{\{  \underset{k\leq b}{\min}{V}(w_k)+y_1\geq 0,\, V(\omega_b)+y_1\leq 2L\}} - \1_{\{ \underset{k\leq b}{\min}{V}(w_k)  +y_2\geq 0 V(\omega_k)+y_2\leq  2L\}}|\right]
\\
+1-\ee^{-\epsilon}+\overset{d}{\underset{j=1}{\sum}} \E_\Q\big(\ee^{- \sum_{|z|=b}\theta_i\ee^{-\beta_i[V(z)+y_1-L]}\1_{\{ V(z)+y_1\leq 2L \}}}-\ee^{ - \sum_{|z|=b}\theta_i\ee^{-\beta_i[V(z)+y_2-L]}\1_{\{ V(z)+y_2\leq 2L \}}}\big).
\end{eqnarray*}
Then we easily deduce that for any ${\bm \theta}\in (\r_{+}^*)^d$, $y\mapsto F_{L,b}({\bm \theta},y)$ is Riemann integrable.
\hfill $\Box$

\begin{Lemma}[Aïdékon \cite{Aid11}]
\label{renewal}
Let $(r_n)_{n\geq 0}$ and $(\lambda_n)_{n\geq 0}$ be two sequences of numbers resp. in $\r_+$ and in (0,1) and such that resp. $\lim_{n\to\infty}\frac{r_n}{n^\frac{1}{3}}=0$, and $\frac{2}{5}<\liminf_{n\to\infty}\lambda_n\leq \limsup_{n\to\infty}\lambda_n<1$. Let $F:\r_+\to \r$ be a Riemann integrable function. We suppose that there exists a non-increasing function $\overline{F}:\r_+\to\r$ such that $|F(x)|\leq \overline{F}(x)$ for any $x\geq 0$ and $\int_{x\geq 0}x\overline{F}(x)<\infty$. Then as $n\to \infty$,
\begin{equation}
\E\big[F(S_n-y_2); \underset{k\leq n}{\min}{S}_k\geq -y_1,\underset{k\in [\lambda_n n,n]}{\min}{S}_{k}\geq y_2\big]\sim \frac{C_-C_+\sqrt{\pi}}{\sigma \sqrt{2}}R(y_1)n^{-\frac{3}{2}}\int_{x\geq 0}F(x)R_-(x)dx
\end{equation}
uniformly in $y_1,\, y_2\in [0,r_n]$
\end{Lemma}
\noindent{\it Proof of Lemma \ref{renewal}.} Lemma \ref{renewal} is a simple extension of the Lemma 2.3 in \cite{Aid11}. By the Markov property, observe that
\begin{eqnarray*}
\E\left[F(S_n-y_2); \underset{k\leq n}{\min}{S}_k\geq -y_1,\underset{k\in [\lambda_n n,n]}{\min}{S}_{k}\geq y_2\right]= \sum_{k=0}^{ \lambda_n n}\E\left[\Upsilon_n(y_2-S_k,k); \underset{j\leq k-1}{\min}{S}_{j}>S_k\geq -y_1\right],
\end{eqnarray*}
(when $k\geq \lambda_n n$, $S_{\lambda_n n}\in [-y_1,0)$ is impossible because of $\underset{k\in [\lambda_n n,n]}{\min}{S}_{k}\geq y_2\geq 0 $), with 
\begin{equation}
\Upsilon_n(y,k):= \E\left[F(S_{n-k}-y);  \underset{j\leq n-k}{\min}{S}_{j}  \geq 0,  \underset{j\in [\lambda_n n -k,n-k]}{\min}{S}_{j}\geq y \right],\quad y\in \r,\, k\in \{0,...,\lambda_n n-1\}.
\end{equation}
Let $(m_n)_{n\in \N}$ a sequence of integers such that $\frac{n}{m_n}$ and $\frac{m_n}{r_n^2}$ go to infinity. First we will show that
\begin{equation}
\label{m_nnegligble} \sum_{k=m_n}^{ \lambda_n n-1}\E\left[
\Upsilon_n(y_2-S_k,k); \underset{j\leq k-1}{\min}{S}_{j}  >S_k\geq -y_1\right]=o(R(y_1)n^{-\frac{3}{2}}).
\end{equation}
We divide the proof of (\ref{m_nnegligble}) in two steps:

- According to the proof of Lemma 2.3 in \cite{Aid11}, for any $k\in [m_n,\frac{n}{3}]$, $y\in [0,2r_n]$, we have
$ \Upsilon_n(y,k)\leq \frac{c}{n^\frac{3}{2}}\sum_{j\geq 0}\overline{F}(j)j $. Thus we deduce that
\begin{eqnarray*}
\sum_{k=m_n}^{n/3}\E\left[ \Upsilon_n(y_2-S_k,k); \underset{j\leq k-1}{\min}{S}_{j}  >S_k\geq -y_1\right]&\leq& c\sum_{k=m_n}^{n/3}\E\left[ \underset{j\leq k-1}{\min}{S}_{j} >S_k\geq -y_1\right] \frac{1}{n^\frac{3}{2}} \sum_{j\geq 0}\overline{F}(j)j
\end{eqnarray*}
Recall that $\sum_{j\geq 0}\overline{F}(j)j<\infty $. Moreover by (\ref{2.6}) (after reversing the time)
\begin{eqnarray*}
\sum_{k=m_n}^{n/3}\E\left[  \underset{j\leq k-1}{\min}{S}_{j} >S_k\geq -y_1\right]\leq  c(1+y_1)^2 \sum_{k=m_n}^{n/3} k^{-\frac{3}{2}}&\leq &c(1+y_1)\frac{1+y_1}{\sqrt{m_n}} 
\\
&=& o(R(y_1)),
\end{eqnarray*}
where we have used that $y_1\leq r_n$ and $\frac{m_n}{r_n^2}$ go to infinity. So $ \sum_{k=m_n}^{n/3}\E[ \Upsilon_n(y_2-S_k,k);\underset{j\leq k-1}{\min}{S}_{j}  >S_k\geq -y_1]=o(\frac{R(y_1)}{n^{\frac{3}{2}}} )$.

- For the second step, we notice that for any $k\in [\frac{n}{3}, \lambda_n n-1]$, $y\in [0,2r_n]$, we have
$ \Upsilon_n(y,k)\leq \sum_{j\geq 0}\overline{F}(j) \P(S_{n-k}-y\in [j-1,j],\, \underline{S}_{n-k}\geq 0)\leq c\frac{r_n}{n^\frac{3}{2}}\sum_{j\geq 0}\overline{F}(j)j.$ Thus we deduce that
\begin{eqnarray*}
\sum_{k=n/3}^{\lambda_n n-1}\E\left[ \Upsilon_n(y_2-S_k,k);   \underset{j\leq k-1}{\min}{S}_{j}  >S_k\geq -y_1\right] &\leq& \sum_{k=n/3}^{\lambda_n n-1}  c\P\left[ \underset{j\leq k-1}{\min}{S}_{j}   >S_k\geq -y_1  \right] \frac{r_n}{n^{\frac{3}{2}}} 
\\
&\leq & c'(1+y_1) \frac{r_n}{n^{\frac{1}{2}}} \frac{1+y_1}{n^\frac{3}{2}}=o(\frac{R(y_1)}{n^\frac{3}{2}}).
\end{eqnarray*}
The proof of (\ref{m_nnegligble}) is now completed. Via Lemma 2.3 in \cite{Aid11}: uniformly in $y_1,\, y_2\in [0,r_n]$, 
\begin{eqnarray}
\nonumber  &&\sum_{k=0}^{m_n}\E\left[ \Upsilon_n(y_2-S_k,k); \underset{j\leq k-1}{\min}{S}_{j}   >S_k\geq -y_1\right]
\\
\nonumber &&=  \sum_{k=0}^{m_n} \P\left[\underset{j\leq k-1}{\min}{S}_{j}    >S_k\geq -y_1  \right] [\frac{C_-C_+\sqrt{\pi}}{\sigma \sqrt{2}}n^{-\frac{3}{2}}\int_{x\geq 0}F(x)R_-(x)dx +o(n^{-\frac{3}{2}})]
\\
\label{m_nnonnegligble}&&= \frac{C_-C_+\sqrt{\pi}}{\sigma \sqrt{2}}R(y_1)n^{-\frac{3}{2}}\int_{x\geq 0}F(x)R_-(x)dx +o(\frac{R(y_1)}{n^\frac{3}{2}}).
\end{eqnarray}
By combining \ref{m_nnegligble} and (\ref{m_nnonnegligble}) we obtain Lemma \ref{renewal}.\hfill$\Box$
\\

Now we end the proof of Proposition \ref{Propconclusion2}. Recalling that $R(x)\underset{x\to\infty}{\sim} c_0 x$ and defining $\rho_{{\bm \beta}, L,b}({\bm \theta})=c_0\rho^{*}_{{\bm \beta}, L,b}({\bm \theta})$, via Lemma \ref{Riemann} and \ref{renewal} and the inequality (\ref{petitpoint}),  we have obtained that: for any ${\bm \theta}\in (\r_+^*)^d,\,  \epsilon>0$,  there exists $(B,L_0)$ such that for any $b\geq B,\, L\geq L_0$, there exists $(A,N)_{(\epsilon)}\in \r_+\times \N$  such that $\forall  n>N$ and $ x\in[A,\frac{3}{2}\log n-A]$, we have
\begin{equation}
\label{deuxsansDelt}
\left|\frac{\ee^x}{x}\E\left(1-\exp\{-\sum_{i=1}^d \theta_i \ee^{-\beta_i x} \widetilde{W}_{n,\beta_j}\} \right)-\rho_{{\bm \beta},L,b}({\bm \theta })\right|\leq \epsilon. 
\end{equation}

 In addition by (\ref{equplowbound}) : {\it There exist $c_6^{({\bm \theta})},\, c_7^{({\bm \theta})}>0$ and  $A,\, N>0$ large such that: {\bf  for any $n\geq N,\, x \in [A,\frac{3}{2}\log n-A]$, 
\begin{equation}
\label{utillower}  c_6^{({\bm \theta})}x\ee^{-x}\leq    \E\left(1-\exp\{-\sum_{i=1}^d \theta_i \ee^{-\beta_i x} \widetilde{W}_{n,\beta_j}\} \right)\leq c_7^{({\bm \theta})} x\ee^{-x}.
\end{equation}}}
For any $ p>0$, let $(L,b)_p$ such that (\ref{deuxsansDelt}) is true (in the sense that: there exists $(A,N)_{(\frac{1}{p})}\in \r_+\times \N$ such that $\forall n>N$ and ...) with $\epsilon=\frac{1}{p}$, then we clearly have $\rho_{\bm \beta}^{(p)} ({\bm \theta}):=\rho_{{\bm \beta},(L,b)_p}({\bm \theta})\in [\frac{c_6}{2}^{({\bm \theta})}, 2c_7^{({\bm \theta})}]$ for any $p> \frac{2}{c_6^{(\theta)}} $. Let $\phi:\N\to \N$ strictly increasing such that $\rho_{\bm \beta}^{(\phi(p))}(\bm \theta)\to \rho_{\bm \beta}({\bm \theta})\in  [\frac{c_6}{2}^{({\bm \theta})}, 2c_7^{({\bm \theta})}]$. 

We shall complete the proof of (\ref{sansDelta2}) by using the following observation: in the display (\ref{deuxsansDelt}) the expectation only depends on $x$ and $n$ whereas $\rho_{{\bf \beta}, L,b}({\bf \theta})$ is independent of $x,n$. 

Fix $\epsilon>0$. Let $p_0>0$ such that $| \rho_{\bm \beta}^{(\phi(p_0))}(\bm \theta) - \rho_{\bm \beta}({\bm \theta}) |\leq \epsilon$ and $\frac{1}{\phi(p_0)}\leq \epsilon $. {\bf Then it suffices to choose (according to (\ref{deuxsansDelt}) and associated to $(L,b)_{\Phi(p_0)}$) $(A,N)_{(\frac{1}{\phi(p_0)})}>0$ such that for any $n\geq N,\, x\in [A,\frac{3}{2}\log n-A]$,
\begin{eqnarray*}
&&\left|\frac{\ee^x}{x}\E\left(1-\exp\{-\sum_{i=1}^d \theta_i \ee^{-\beta_i x} \widetilde{W}_{n,\beta_j}\} \right)-\rho_{{\bm \beta}}({\bm \theta })\right|
\\
&&\leq \left|\frac{\ee^x}{x}\E\left(1-\exp\{-\sum_{i=1}^d \theta_i \ee^{-\beta_i x} \widetilde{W}_{n,\beta_j}\} \right)- \rho_{\bm \beta}^{(\phi(p))}(\bm \theta)\right|+  | \rho_{\bm \beta}^{(\phi(p))}(\bm \theta) - \rho_{\bm \beta}({\bm \theta}) |\leq 2\epsilon.
\end{eqnarray*}}
This completes the proof of (\ref{sansDelta2}). We stress here that our argument shows that all the possible extractions of $ \rho_{\bm \beta}^{(p)} ({\bm \theta})$ converge, in the end, to the same limit. 
\\

To complete the proof of the Proposition \ref{Propconclusion2} it remains to prove that $\rho_{\bm \beta}$ is a continuous function at $0$. For any ${\bm \theta}\in \r_{+,*}^d$ let $i*\in \{1,...,d\}$ such that $\theta_{i*}:=\underset{i\in [1,d]}{\max} \theta_i^\frac{1}{\beta_i}$. By applying $d$ times the inequality $1-\ee^{x+y} \leq1-\ee^x+1-\ee^y$ with $x,\, y\geq 0$, then Corollary \ref{partieneg},  we deduce that there exists $A,N>0$ such that for any $n\geq N$ and $x\in [A+ \frac{1}{\beta_{i*}}\log \theta_{i*},\frac{3}{2}\log n-A+ \frac{1}{\beta_{i*}}\log \theta_{i*}]$,
\begin{eqnarray*}
&&\frac{\ee^x}{x}\E\left(1-\exp\{-\sum_{i=1}^d \theta_i \ee^{-\beta_i x} \widetilde{W}_{n,\beta_j}\} \right)\leq  \sum_{i=1}^d \frac{\ee^x}{x}\E\left(1-\exp\{- \ee^{-\beta_i (x-\frac{1}{\beta_{i}}\log \theta_{i})} \widetilde{W}_{n,\beta_{i}}\} \right)
\\
&&\leq  \sum_{i=1}^d\frac{x-\frac{1}{\beta_{i}}\log \theta_{i}}{x} {\ee^{ \frac{1}{\beta_{i}}\log \theta_{i}} }       \frac{\ee^{x- \frac{1}{\beta_{i}}\log \theta_{i}}}{x-\frac{1}{\beta_{i}}\log \theta_{i}}\E\left(1-\exp\{- \ee^{-\beta_i (x-\frac{1}{\beta_{i}}\log \theta_{i})} \widetilde{W}_{n,\beta_{i}}\} \right)
\\
&&\leq 2c_2 d\theta_{i*}^\frac{1}{\beta_{i}}\to 0,\quad \text{when } \theta_{i*}\to 0.
\end{eqnarray*}
It implies that $\rho_{\bm \beta}({\bm \theta}) \underset{{\bm \theta}\to 0}{\to }0$, thus the proof of Proposition of \ref{Propconclusion2} is terminated.
\hfill$\Box$

\appendix

\appendix

\section{Auxiliary estimates}
Here we prove a slight extension of the Lemma B.3 in \cite{Aid11}. It will be used to prove the Lemma \ref{ExtensAid1}.

For $\alpha>0$, $d_1,\, d_2,\, a\geq 0$, $n\geq 1$ and $0\leq i\leq n$, we define
\begin{equation}\label{defk}
k_i^{(d_1,d_2)}(x,a)=k_i^{(d_1,d_2)} :=
\begin{cases}
-d_1+i^\alpha , &\hbox{if $0\le i\le \lfloor \log d_1\rfloor$,} \cr
i^{\alpha}-x, &\hbox{if $ \lfloor \log d_1 \rfloor \le i\le \lfloor n/2\rfloor$,} \cr
a+(n-i)^{\alpha}-d_2, &\hbox{if $\lfloor n/2 \rfloor <i\le n$.}\cr
\end{cases}	
\end{equation}

\begin{Lemma}\label{A.111}
Let $\alpha \in (0,1/6)$ and $\epsilon>0$. There exist $d_1,\, d_2>0$ large enough such that for any $u,\, x\ge 0$, $a\in (0,10\log n)$ and $n\geq \ee^{d_1+d_2}$,
\begin{eqnarray}\label{eq:a6}
&&\P\Big\{\exists   i\le   n :\, S_i\le  k_i^{(d_1,d_2)},\, \min_{j\le n} S_j \ge -x, \;
 \min_{\lfloor  n/2\rfloor < j\le n} S_j \ge a, \;S_n \le a +u\Big\}\\
\nonumber &\le& (1+x)(1+u)^2 {\epsilon\over n^{3/2}} .
\end{eqnarray}
\end{Lemma}
\noindent {\it Proof}. 
We treat $n/2$ as an integer. Let $E$ be the event in (\ref{eq:a6}). We have
$\P(E)\le \sum_{i=1}^n \P(E_i)$ where
$$
E_i :=\{  S_i\le  k_i^{(d_1,d_2)},\, \min_{j\le n} S_{j} \ge  -x, \;
 \min_{ \frac{n}{2} n <j\le n} S_j \ge a, \;S_n \le a +u\}.
$$

When $d$ is large enough, by the Markov property at time $ i\in[ 1, \lfloor \log d_1 \rfloor]$ and (\ref{2.7}),
\begin{eqnarray}
\nonumber \sum_{i=1}^{\lfloor \log d_1 \rfloor} \P(E_i)&\leq& \sum_{i=1}^{\lfloor\log d_1\rfloor}  {c_{54}(1+u)^2\over n^{3/2}}\E\left[ (1+S_i+x)_+{\bf 1}_{\{  S_i\le  -\frac{d_1}{2}\}}\right]
\\
\nonumber &\leq& \sum_{i= 0  }^{\lfloor\log d_1\rfloor}  {c_{54}(1+u)^2\over n^{3/2}}(1+x) \P(S_i\leq - \frac{d_1}{2}) 
\\
 &\leq& \sum_{i=0}^{\lfloor\log d_1\rfloor}  {c_{54}(1+u)^2\over n^{3/2}}(1+x)4\sigma^2 i/d_1 \leq \epsilon (1+x)(1+u)^2n^{-\frac{3}{2}}.
\end{eqnarray}

\noindent Now we treat the case $i\in ( \lfloor\log d_1\rfloor,n/2]$. By the Markov property at time $i\ge 1$,  (\ref{2.6}) and (\ref{2.7}), we have
\begin{equation}
 \P(E_i) \le \left\{ \begin{array}{ll}  {c_{21}(1+u)^2\over n^{3/2}}\E\left[ (1+S_i+x)_+{\bf 1}_{\{  S_i\le  i^{\alpha} -x,\, \min_{j\le i} S_j \ge  -x\}}\right]\,\, &\text{if}\quad \lfloor\log d_1\rfloor \leq i\leq n/3,
\\
{c_{12}(1+a+u)^2\over n^{3/2}}\E\left[ (1+S_i+x)_+{\bf 1}_{\{  S_i\le  i^{\alpha} -x,\, \min_{j\le i} S_j \ge  -x\}}\right]\,&\text{if}\quad  n/3\leq i\leq \frac{n}{2}.
\end{array}\right. 
\end{equation}
By (\ref{2.7}), recalling that $a\leq \log n$, it yields that for $d_1$ large enough,
\begin{equation}\label{eq:a6sum1}
\sum_{i= \lfloor\log d_1 \rfloor }^{n/2} \P(E_i) \leq (1+x)[ \sum_{i=  \lfloor\log d_1 \rfloor}^{n/3} {c_{22}(1+u)^2\over n^{3/2}} \, \,{ (1+i^{\alpha})^3\over i^{3/2}}+ \sum_{i=n/3}^{n/2} {c_{23}(1+u+a)^2\over n^{3/2} n^{3/2-3\alpha}}  ]\leq  (1+u)^2(1+x){\varepsilon \over n^{3/2}},
\end{equation}
Finally we treat the case $i\in [n/2,n]$. By the Markov property at time $i$ and (\ref{2.6}), we have
$$
\P(E_i) \le {c_{24}(1+u)^2\over (n-i+1)^{3/2}}\E\left[ (1+S_i-a)_+{\bf 1}_{\{  S_i\le  a+(n-i)^{\alpha}-d_2,\, \min_{j\le n} S_j \ge-x,\,\min_{ \frac{n}{2} n <j\le i} S_j \ge a \}}\right].
$$
If $i\geq n-d_2^\frac{1}{\alpha}$, clearly we have $\P(E_i)=0$. If $n-d_2^\frac{1}{\alpha}\ge i\ge  2n/3$, we use (\ref{2.6}) to see that $ \P(E_i) \le c_{25}(1+x)(1+u)^2{(1+n-i)^{3\alpha-{3\over 2}}\over n^{3/2}}$. Therefore for $d_2$ large enough we have
\begin{equation}\label{eq:a6sum2}
\sum_{i=\lfloor 2n/3 \rfloor}^{n} \P(E_i) \le \sum_{i=\lfloor 2n/3 \rfloor}^{n-d_2^\frac{1}{\alpha}} c_{25}(1+x)(1+u)^2{(1+n-i)^{3\alpha-{3\over 2}}\over n^{3/2}}  \leq  (1+x)(1+u)^2{\varepsilon \over n^{3/2}}.
\end{equation}
If $n/2< i< 2 n/3$, we simply write
\begin{eqnarray*}
\P(E_i)
&\le&
{c_{26}(1+u)^2\over (n-i+1)^{3/2}}\E\left[ (1+S_i-a){\bf 1}_{\{  a\le S_i\le  a+(n-i)^{\alpha},\, \min_{j\le i} S_j \ge  -x\}}\right]\\
&\le& c_{27} (1+u)^2 {(n-i)^{\alpha} \over (n-i+1)^{3/2}} \P(a\le S_i\le  a+(n-i)^{\alpha},\, \min_{j\le i} S_j \ge  -x)\\
&\le&
c_{28} (1+x)(1+u)^2 {n^{\alpha}(a+n^{\alpha})^2\over n^3}.
\end{eqnarray*}
by (\ref{2.6}). We deduce that
\begin{equation}\label{eq:a6sum3}
\sum_{i=n/2}^{\lfloor 2n/3\rfloor} \P(E_i) \le c_{29}(1+x)(1+u)^2{ (n^{\alpha} +a)^2 \over n^{2-\alpha}}.
\end{equation}
Inequalities (\ref{eq:a6sum1}), (\ref{eq:a6sum2}) and (\ref{eq:a6sum3}) suffice to prove the Lemma \ref{A.111}
\hfill $\Box$

\section{Good vertex}
The Lemma \ref{ExtensAid1} below gives a control on the behavior of the different brothers of the spine. It will be used to prove the Lemma \ref{t17} of the Appendix C. It allwo

Let us recall some notations. Let $(e_k,0\leq k \leq n)$ such that
\begin{equation}
\label{defe_k}
e_k=e_k^{(n)}=\left\{
         \begin{array}{ll}
            k^{\frac{1}{12}}& \qquad \mathrm{if}\quad 0<k\leq \frac{1}{2}n, \\
            (n-k)^{ \frac{1}{12}}&\qquad \mathrm{if} \quad \frac{n}{2}<k\leq n, \\
         \end{array}
       \right.
\end{equation}
and denote
\begin{equation}
\label{defd_k}
d_k=d_k(n,x,L,B_1):=\left\{
         \begin{array}{ll}
         -B_1 &\qquad \text{if}\quad 0\leq k\leq \log B_1, \\
           -x  & \qquad \mathrm{if}\quad \log B_1\leq k\leq\frac{n}{2}, \\
           \max(a_n-x-L-1,0) &\qquad \mathrm{if} \quad \frac{n}{2}<k\leq n. \\
         \end{array}
       \right..
\end{equation}
We say that $|u|=n$ is a $(x,L,B_1,B_2)-$good vertex if $u\in \mathcal{Z}^{x,x+L,x-L}$ and

\begin{equation}
\label{defgood_k}
\underset{w\in\Omega(u_k)}{\sum}\ee^{-(V(v)-d_k)}\{1+(V(v)-d_k)_+\}\leq B_2\ee^{-e_k},\quad \quad \forall 1\leq k\leq n.
\end{equation}

\begin{Lemma}[Aïdékon \cite{Aid11}]
\label{ExtensAid1}
Fix $L\geq 0$. For any $\epsilon >0$, we can find $B_1^{(0)},B_2^{(0)} $ large enough such that for any $B_1\geq B_1^{(0)} ,\, B_2\geq B_2^{(0)}$, as in (\ref{defe_k}) and (\ref{defd_k}), for any $n\geq \ee^{B_1+B_2} $ and $x\geq 0$
\begin{equation}
\label{toprove}
\Q\left(w_n\, \text{ is not a } (x,L,B_1,B_2)\text{-good vertex},\, w_n\in \mathcal{Z}_n^{x,x+L,x-L}\right) \leq (1+x)\epsilon n^{-\frac{3}{2}}.
\end{equation}
\end{Lemma}
\noindent{\it Proof of Lemma \ref{ExtensAid1}.}
From Lemma \ref{A.111} there exists $B(L),\, c(L)>0$ large enough such that for any $B_1\geq B$, $n\geq  \ee^{B_1+B_2}$, and $x\geq 0$,
\begin{equation}
\label{autoi}
\Q\Big( \{w_n\in \mathcal{Z}_n^{x,x+L,x-L}\}\cup \overset{\frac{n}{2}}{\underset{j=0}{\cup}} \{V(w_j)\leq d_{j}+2e_{j}\} \underset{j=\frac{n}{2}+1}{\overset{n}{\cup}}  \{V(w_j)\leq d_{j}+2e_{j}-c(L)\}\Big) \leq \frac{\epsilon}{n^\frac{3}{2}}(1+x),
\end{equation}
(with according to the notation of Lemma \ref{A.111}, $a\leftrightarrow  a_n-x-L-1$, $u \leftrightarrow 2L$, $x \leftrightarrow x$, $B_1 \leftrightarrow d_1$ , $c(L)\leftrightarrow  d_2$,  $i^\alpha \leftrightarrow  2e_i$). From now we fix $c(L)>0$. For any $n,\, B_1>0$, we denote $\mathcal{H}_n $ the event $( \overset{\frac{n}{2}}{\underset{j=0}{\cup}} \{V(w_j)\leq d_{j+1}+2e_{j+1})\} \underset{j=\frac{n}{2}+1}{\overset{n}{\cup}}  \{V(w_j)\leq d_{j+1}+2e_{j+1}-c(L))\})^c $. Consequently, it is enough to show that for $B_1,\, B_2$ large enough,
\begin{eqnarray}
\nonumber \sum_{k=1}^n \Q\Big(\sum_{v\in \Omega(w_k)}\ee^{-(V(v)-d_k)}\{1+(V(v)-d_k)_+\}>B_2\ee^{-c(L)/2}\ee^{-\frac{V(w_{k-1})-d_k}{2}}, w_n\in \mathcal{Z}_n^{x,x+L,x-L},\,  \mathcal{H}_n \Big)
\\
\label{C.4}\leq \epsilon(1+x) n^{-\frac{3}{2}}.
\end{eqnarray}
We see that
\begin{eqnarray*}
&&\sum_{v\in \Omega(w_k)}\ee^{-(V(v)-d_k)}(1+(V(v)-d_k)_+)
\\
&&\leq \ee^{-(V(w_{k-1})-d_k)} \sum_{v\in \Omega(w_k)}\ee^{-(V(v)-V(w_{k-1}))}\{ 1+(V(w_{k-1})-d_k)_+ + (V(v)-V(w_{k-1}))_+\}
\\
&&\leq \ee^{-(V(w_{k-1})-d_k)}(1+(V(w_{k-1})-d_k)_+)\sum_{v\in \Omega(w_k)} \ee^{-(V(v)-V(w_{k-1}))}\{ 1+ (V(v)-V(w_{k-1}))_+\}.
\end{eqnarray*}
By denoting for any $|v|\geq 1$, $\xi(v):= \sum_{w\in \Omega(v)}(1+(V(w)-V(\overset{ \leftarrow}{v}))_+)\ee^{-V(w)-V(\overset{ \leftarrow}{v}))}$, we have then
$$\sum_{v\in \Omega(w_k)}\ee^{-(V(v)-d_k)}(1+(V(v)-d_k)_+)\leq \ee^{(V(w_{k-1})-d_k)}(1+(V(w_{k-1})-d_k)_+)\xi(w_k). $$
Equation (\ref{C.4}) boils down to showing that , for $B_1,\, B_2$ large enough,
\begin{eqnarray}
\nonumber&&\sum_{k=1}^n\Q\Big(\xi(w_k)> B_2 \ee^{-c(L)/2}\frac{\ee^{\frac{V(w_{k-1})-d_k}{2} }}{1+(V(w_{k-1})-d_k)_+}, w_n\in \mathcal{Z}_n^{x,x+L,x-L},\,  \mathcal{H}_n \Big)\leq 
\\
\label{C.5}
&&\sum_{k=1}^n\Q\Big(\xi(w_k)> B_2 \ee^{-c(L)/2}\ee^{\frac{V(w_{k-1})-d_k}{3} }, w_n\in \mathcal{Z}_n^{x,x+L,x-L},\,   \mathcal{H}_n \Big)\leq \epsilon(1+x)n^{-\frac{3}{2}} .
\end{eqnarray}

First we deal with the case $k\in [1,\frac{3n}{4}]$. By the Markov property at time $k$, we get 
\begin{eqnarray*}
\Q\Big(\xi(w_k)> B_2\ee^{-c(L)/2} \ee^{\frac{V(w_{k-1})-d_k}{3} }, w_n\in \mathcal{Z}_n^{x,x+L,x-L},\,  \mathcal{H}_n \Big) 
\\
 \leq   \Q\Big[ \lambda(V(w_k),k,n)\1_{\{ \xi(w_k)>  B_2\ee^{-c(L)/2} \ee^{\frac{V(w_{k-1})-d_k}{3} },\,  V(w_j)\geq -x,\forall j\leq k \}} \Big],
\end{eqnarray*}
where 
\begin{eqnarray*}
 \lambda(r,k,n):= \P_{r}(S_j\geq d_{j+k}+e_{j+k},\, \forall j\leq \frac{n}{2}-k,\,  S_j\geq d_{j+k}+e_{j+k}-c(L),\, \forall j\in ( \frac{n}{2}-k, n-k],\, \qquad \\ S_{n-k} \leq a_n-x+L,
 \min_{j\in [k,\frac{n}{2}]}S_{j-k}\geq -x,\, \min_{j\in [n/2,n]}S_{j-k}\geq a_n-x-L-1).
\end{eqnarray*}
When $k\in [1,n/3]$ by (\ref{2.7}),
\begin{eqnarray}
\nonumber \lambda(r,k,n)&\leq&  \P_{r}( S_{n-k} \leq a_n-x+L, 
 \min_{j\in [k,\frac{n}{2}]}S_{j-k}\geq -x,\, \min_{j\in [n/2,n]}S_{j-k}\geq a_n-x-L-1)
\\
\label{Cfouca0}  &\leq&c_{30}(1+L)^2n^{-\frac{3}{2}}(1+(r+x)_+).
\end{eqnarray}
When $k\in (n/3,3n/4]$ by (\ref{2.6}) (recalling $d_{j+k}+e_{j+k}-c(L)\geq a_n-x-L-1$ for any $k\in (n/3,3n/4],\, j\in [0,n-k-\frac{n}{5}]$),
\begin{eqnarray}
\nonumber  \lambda(r,k,n)&\leq&  \P_{r}( S_{n-k} \leq a_n-x+L,  \min_{j\in [k,n]}S_{j-k}\geq a_n-x-L-1) 
\\
\label{Cfouca} &\leq& c_{30}(1+L)^2n^{-\frac{3}{2}}(1+(r+x)_+).
\end{eqnarray}
This yields that
\begin{eqnarray}
\nonumber&&\Q\Big(\xi(w_k)>  B_2 \ee^{-c(L)/2}\ee^{\frac{V(w_{k-1})-d_k}{3} }, w_n\in \mathcal{Z}_n^{x,x+L,x-L},\, \mathcal{H}_n \Big)
\\
 \nonumber &&\leq    c_{30}(1+L)^2n^{-\frac{3}{2}} \Q\Big[(1+(V(w_k)+x)_+) \1_{\{ \xi(w_k)> B_2 \ee^{-c(L)/2}\ee^{\frac{V(w_{k-1})-d_k}{3} },\,  V(w_j)\geq -x+(e_j-c(L))_+,\forall j\leq k \}}  \Big]
\\
\nonumber &&\leq    c_{30}(1+L)^2n^{-\frac{3}{2}} \Q_x\Big[(1+V(w_k)_+) \times
\\
\label{C.6} &&\qquad\qquad\qquad  \1_{\{ \xi(w_k)> B_2 \ee^{\frac{1}{4}(V(w_{k-1})+(B_1-x)\1_{\{ k-1\leq \log B_1\}}) -c(L)/2 },\,  V(w_j)\geq0,\forall j\leq k \}}  \Big],
\end{eqnarray}
indeed, $V(w_k)+x\geq (e_k-c(L))_+$ implies $\frac{V(w_k)-d_k}{3} \geq \frac{V(w_k)+x}{4}$ $\forall k\geq \log B_1$.  On the other hand, we have
$$ 1+(V(w_k)+x)_+\leq 1+(V(w_{k-1})+x)_+ + (V(w_k)-V(w_{k-1}))_+. $$
Let $(\xi,\Delta)$ be generic random variable distributed as $(\xi(w_1),V(w_1)_+)$ under $\Q$ and independent of the other random variables. By the Markov property at time $k-1$, we obtain that
\begin{eqnarray*}
\Q_x\Big[(1+V(w_k)_+) \1_{\{ \xi(w_k)> B_2 \ee^{\frac{1}{4}(V(w_{k-1})+(B_1-x)\1_{\{ k-1\leq \log B_1\}}) -c(L)/2 }  ,\,  V(w_j)\geq 0,\forall j\leq k \}}  \Big]
\\
\leq \Q_x\Big[\kappa_k(V(w_{k-1}))\1_{\{ V(w_j)\geq 0,\, \forall j\leq k-1\}}\Big],
\end{eqnarray*}  
with, for any $z\geq 0$
$$ \kappa_k(z):= \left\{ \begin{array}{ll} (1+z) \1_{\{ \xi>   B_2\ee^{-c(L)/2} \ee^{\frac{z-x+B_1}{4} }  \}} +\Delta_+\1_{\{ \xi>  B_2\ee^{-c(L)/2} \ee^{\frac{z-x+B_1}{4} }  \}}&\text{if } k\leq \log B_1
\\
(1+z) \1_{\{ \xi>   B_2\ee^{-c(L)/2} \ee^{\frac{z}{4} }  \}} +\Delta_+\1_{\{ \xi>  B_2\ee^{-c(L)/2} \ee^{\frac{z}{4} }  \}}&\text{if } k> \log B_1
\end{array}\right. .
$$
 In view of (\ref{C.6}), it follows that 
\begin{equation}
\label{thomztio}
\sum_{k=1}^{ \frac{3n}{4}} \Q\Big(\xi(w_k)> B_2 \ee^{-c(L)/2} \ee^{\frac{V(w_{k-1})-d_k}{3} }, w_n\in \mathcal{Z}_n^{x,x+L,x-L},\mathcal{H}_n \Big) \leq    c_{30}(1+L)^2n^{-\frac{3}{2}}(D_1+D_2),
\end{equation}
where $c_{30}$ is the constant defined in (\ref{Cfouca0}) and (\ref{Cfouca}) and
\begin{eqnarray*}
&&D_1:=  \sum_{k=1}^{ \log B_1} \Q_x\Big[ (1+V(w_k))\1_{\{V(w_k) +B_1-x\leq 4(\log \xi- \log B_2) +2c(L) \}},\, \underset{j\leq k}{\min}{V}(w_j)\geq 0\Big]+ 
\\
&&\qquad\qquad\qquad\qquad\qquad  \sum_{k=\log B_1+1 }^{ \frac{3n}{4}} \Q_x\Big[ (1+V(w_k))\1_{\{ V(w_k)\leq 4(\log \xi- \log B_2) +2c(L) \}},\, \underset{j\leq k}{\min}{V}(w_j)\geq 0\Big],
\\
&&D_2:=\sum_{k=1}^{  \log B_1} \Q_x\Big[ \Delta_+\1_{\{ V(w_k)+B_1-x\leq 4(\log \xi- \log B_2) +2c(L) \}},\,  \underset{j\leq k}{\min}{V}(w_j)\geq 0\Big]+
\\
&&\qquad\qquad\qquad\qquad\qquad \sum_{k=\log B_1+ 1}^{ \frac{3n}{4}} \Q_x\Big[ \Delta_+\1_{\{ V(w_k)\leq 4(\log \xi- \log B_2) +2c(L) \}},\,  \underset{j\leq k}{\min}{V}(w_j)\geq 0\Big]
\end{eqnarray*}
When $k\in [1,\log B_1]$, using the definition of $d_k$ in (\ref{defd_k}), observe that (for $B_2\geq \ee^{10 c(L)}$)
\begin{eqnarray}
\nonumber &&\Q_x\Big[ (1+V(w_k))\1_{\{ V(w_k)+B_1-x\leq 4(\log \xi- \log B_2)+ 2c(L) \}},\,  \underset{j\leq k}{\min}{V}(w_j)\geq 0 \Big]
\\
\nonumber && \leq \E_\Q\Big[ (1+4\log \xi+x)\1_{\{B_1 \leq 3\log \xi- V(w_k) \}}\Big]
\\
\label{over1}&& \leq c_{30} (1+x)\Q(B_1\leq 4\log \xi-V(w_k))+ \E_{\Q}(\log \xi\1_{\{B_1\leq 4\log \xi-V(w_k) \}}),
\end{eqnarray}
and 
\begin{equation}
\label{bibioubiou} \Q\Big[ \Delta_+\1_{\{ V(w_k)+B_1\leq 4(\log \xi- \log B_2)+ 2c(L) \}},\,  \underset{j\leq k}{\min}{V}(w_j)\geq -x\Big] \leq 16\Q\Big[\Delta_+ (4\log(\xi)_+ +\E(|V(w_k)|)\Big]/B_1.
\end{equation}
Recalling that  for any $k\in [1,\log B_1]$, $\E(V(w_k)^2)\leq \sigma^2 \log B_1$, and $\E_\Q((1+\log_+ \xi)^2)<\infty$, we deduce that for $B_1(c_{30}(1+L)^2)$ large enough we have
\begin{equation}
   \sum_{k=1}^{ \log B_1} \Q_x\Big[ (1+V(w_k))\1_{\{V(w_k) +B_1-x\leq 4(\log \xi- \log B_2) +2c(L) \}},\, \underset{j\leq k}{\min}{V}(w_j)\geq 0\Big] 
\label{ragnignia} \leq \frac{\epsilon(1+x)}{c_{30}(1+L)^2}n^{-\frac{3}{2}},
\end{equation}
 \begin{equation}
\label{ragnignia2}  \sum_{k=1}^{  \log B_1} \Q_x\Big[ \Delta_+\1_{\{ V(w_k)+B_1-x\leq 4(\log \xi- \log B_2) +2c(L) \}},\,  \underset{j\leq k}{\min}{V}(w_j)\geq 0\Big] \leq \frac{\epsilon(1+x)}{c_{30}(1+L)^2}n^{-\frac{3}{2}}.
\end{equation}

\noindent When $k\in (\log B_1, 3n/4]$,
\begin{eqnarray}
\nonumber \sum_{k= \log B_1}^{ 3n/4}\E_{\Q_x}\Big[ (1+V(w_k))\1_{\{ V(w_k)-d_k\leq '(\log \xi- \log B_2)+2c(L)\}},\,  \underset{j\leq k}{\min}{V}(w_j)\geq -x\Big]
\\
\label{interdem} \leq  \sum_{k= 0}^{+\infty}\E_{\Q_x}\Big[ (1+V(w_k))\1_{\{ V(w_k)\leq 4(\log \xi- \log B_2)+2c(L)  \}},\,    \underset{j\leq k}{\min}{V}(w_j)\geq 0\Big] .
\end{eqnarray}
 Notice that in (\ref{interdem}), the term  inside the expectation is $0$ if $\log B_2>  \log \xi+3c(L)/2 $. Therefore, we can add the indicator that $ \log B_2 -3c(L)/2\leq  \log \xi $. By Lemma B.2 (i), we get that
\begin{eqnarray}
\nonumber &&\sum_{k= 0}^{+\infty}\E_{\Q_x}\Big[ (1+V(w_k))\1_{\{ V(w_k)\leq 4(\log \xi- \log B_2)+2c(L)  \}},\, \underset{j\leq k}{\min}{V}(w_j)   \geq 0\Big]  
\\
\nonumber &&\leq c_{30}\Q\big[\1_{\{  \log B_2 -2c(L)\leq  4\log \xi  \}} (1+(\log \xi -\log B_2)_+^2)  \big](1+x)
\\
\label{hatcoum} &&\leq c_{30}\Q\big[\1_{\{  \log B_2 -2c(L)\leq  4\log \xi \}}(1+\log_+\xi)^2 \big](1+x).
\end{eqnarray}
Observe that $\xi\leq X+\tilde{X}$ with the notation of  (\ref{intcdtn}). Going back to the measure $\Q$, we get 
\begin{eqnarray}
\nonumber &&\sum_{k= 0}^{+\infty}\E_{\Q_x}\Big[ (1+V(w_k))\1_{\{ V(w_k)\leq 4(\log \xi- \log B_2)+2c(L)  \}},\, \underset{j\leq k}{\min}{V}(w_j)  \geq 0\Big]
\\
\label{over4} &&\leq \E[X\1_{\{   \log B_2 -2c(L) \leq 4\log(X+\tilde{X})\}}(1+\log_+(X+\tilde{X}))^2 \leq \frac{(1+x)\epsilon}{c_{30}(1+L)^2},
\end{eqnarray}
for $B_2$ large enough. Similarly,
\begin{equation}
\label{over3}\sum_{k=\log B_1}^{ 3n/4} \Q_x[\Delta_+\1_{\{ V(w_k) \leq 4(\log (\xi-\log B_2)+2c(L) \}}]\leq \frac{\epsilon(1+x)}{c_{30}(1+L)^2}.
\end{equation}
In order to prove (\ref{C.4}), by combining (\ref{C.5}), (\ref{ragnignia}), (\ref{ragnignia2}), (\ref{over4}) and (\ref{over3}) it remains to treat the case $\frac{3n}{4}\leq k\leq n$.  Recalling (\ref{C.5}), we want to show that for $B(=B_2\ee^{-c(L)/2})$ large enough, $n\geq \ee^{B}$ and $x\geq 1$, 
\begin{equation}
\label{C.52}
\sum_{k=\frac{3n}{4}}^n\Q\Big(\xi(w_k)> B \ee^{\frac{V(w_{k-1})-d_k}{3} }, w_n\in \mathcal{Z}_n^{x,x+L,x-L},\, \mathcal{H}_n \Big)\leq \epsilon(1+x)n^{-\frac{3}{2}} .
\end{equation}
This case is quasi identical to the proof of (C.8) (in the proof of Lemma C.1 \cite{Aid11}) in \cite{Aid11}, then we omit the details of the proof of (\ref{C.52}).

\hfill$\Box$

\section{Proof of Lemma \ref{t17}.}
\noindent{\textbf{Lemma \ref{t17}}}
{\it $\forall \eta,L>0$ $\exists D(L,\eta)>0$  and  $B(L,\eta) \geq 1$ such that $\forall b\geq B,\, n\geq \ee^{5b},\,x\geq D$,
\begin{equation}
\Q((\xi_n(x,b,L))^c,\omega_n\in \mathcal{Z}_n^{x,x+L,x-L})\leq \eta n^{-\frac{3}{2}} (1+x).
\end{equation}}
\noindent {\it Proof of lemma \ref{t17}.} Let $L,\eta>0$. 
According to Lemma \ref{ExtensAid1}, there exists $B_0$($=B_0(L,\eta)$) such that, for any $B_1,\, B_2\geq B_0$, $n\geq \ee^{B_1+B_2}$ and $x\geq 0$
\begin{equation}
\Q(w_n\in \mathcal{Z}_n^{x,x+L,x-L},w_n\text{ is not a $(x,L,B_1,B_2)$-good vertex })\leq \frac{\eta}{n^{\frac{3}{2}}}(1+x).
\label{t13}
\end{equation}
Now we fix $B_2\geq B_0$. For $\xi_n$ to happen, every brother of the spine at generation less than $n-b$ must have its descendants at time $n$ greater than $a_n(x)+L$. In others words,
\begin{equation}
\Q((\xi_n)^c,\omega_n\text{ is a good vertex})=\Q\Big[1-\overset{n-b}{\underset{k=1}{\prod}}\underset{u\in\Omega(\omega_k)}{\prod}p(u,x-L),\omega_n\text{ is a good vertex}\Big],
\label{t12}
\end{equation}
where $p(u,x-L)=\P(M_{n-|u|}\geq a_n-x-V(u)+L))$ is the probability that the branching random walk rooted at $u$ has its minimum greater $a_n-x+L$ at time $n-|u|$. From Proposition \ref{tenstionaid} \cite{Aid11}, we see that
$$1-p(u,x-L)\leq c_{2}(1+(x+V(u)-L)_+)\ee^{-(x+V(u)-L)},\quad |u|\leq \frac{n}{2}.$$
Moreover as $w_n$ is a good vertex, we have
\begin{eqnarray*}
 \sum_{u\in \Omega(w_k)} (1+(V(u)+x-L)_+)\ee^{-(x-L)-V(u)}\leq \left\{ \begin{array}{ll}
 B_2\ee^{B_1}\ee^{-k^\frac{1}{12}} (1+x) \ee^{-(x-L)} ,\qquad \text{if }\quad k\leq \log B_1, \\
B_2 \ee^{-k^\frac{1}{12}+L} ,\qquad \text{if }\quad k\in (\log B_1,n/2],
\end{array}\right.
\end{eqnarray*}
where we have used that $d_k=-B_1$ when $k\leq \log B_1$  and $d_k= -x$ when $k\in (\log B_1,n/2]$. Using the inequality $p\geq \ee^{\frac{p-1}{2}}$ for $p$ close enough to $1$, it implies that for $x$ large enough and $1\leq k\leq n/2$,
\begin{eqnarray*}
\overset{ n/2}{\underset{k=1}{\prod}}\underset{u\in\Omega(\omega_k)}{\prod}p(u,x-L) &\geq&  \exp(-(2 c_{2} B_2 \ee^{L})(\ee^{B_1} (1+x)\ee^{-x}\overset{\log B_1}{\underset{k=1}{\sum}}\ee^{-k^\frac{1}{12}} + \overset{n/2}{\underset{k=\log B_1}{\sum}}\ee^{-k^\frac{1}{12}}))
\\
&\geq& \exp(- (2 c_{2} B_2 \ee^{L})(c_{31}\ee^{B_1}(1+x)\ee^{-x}  + \ee^{-(\log B_1)^\frac{1}{13}})).
\end{eqnarray*}
Therefore there exists $B_1>0$ and  $D_1(L,B_1)>0$ large enough such that for any $x\geq D_1$ 
\begin{equation}
\overset{ n/2}{\underset{k=1}{\prod}}\underset{u\in\Omega(\omega_k)}{\prod}p(u,x-L)\geq (1-\frac{\eta}{L^2})^{1/2}.
\label{t10}
\end{equation}
From now $B_1$ and $D_1$ are fixed.
\\

If $k>n/2$, since $W_n$ (defined in (\ref{defW})) is a martingale, we have $1=\E[W_l]\geq \E[\ee^{-M_l}]\geq \ee^{-x}\P(M_l\leq x)$ for any $l\geq 1$ and $x\in \r$. We get that
$$1-p(u,x-L)\leq \P\left(M_{n-|u|}\leq a_n-x+L-V(u)\right)\leq \ee^{a_n-x+L}\ee^{-V(u)}.$$
We rewrite it (we have $x-L\geq0$), $1-p(u,x-L)\leq n^{\frac{3}{2}}e^{-V(u)}\ee^{-x+L}=\ee^{-(V(u)-d_k)}\ee^{L}$ for $n/2<k\leq n$. Since $w_n$ is a good vertex, we get that $\underset{u\in\Omega(w_k)}{\prod}p(u,x-A)\geq \ee^{-c_{32}(B_1,B_2)e_k\ee^{2L}}=\ee^{- c_{32}(B_1,B_2)(n-k)^{1/12}\ee^{2L}}$. Consequently,
$$\overset{n-b}{\underset{k=\lfloor n/2 \rfloor+1}{\prod}}\underset{u\in\Omega(\omega_k)}{\prod}p(u,x-L)\geq \ee^{- c_{32}(B_1,B_2)\ee^{L}\overset{n-b}{\underset{\lfloor n/2 \rfloor +1}{\sum}}\ee^{-(n-k)^\frac{1}{12}}}.$$
It yields that there exists $B(L,\eta, c_{32}(B_1,B_2))\geq B_0$ large enough such that $\forall b\geq B$, $n>b$, we have,
\begin{equation}
\overset{n-b}{\underset{k=\lfloor n/2 \rfloor +1}{\prod}}\underset{u\in\Omega(\omega_k)}{\prod}p(u,x-L)\geq (1- \frac{\eta}{L^2})^\frac{1}{2}.
\label{t11}
\end{equation}
In view of (\ref{t10}) and (\ref{t11}), we have for $b\geq B$, $x\geq D_1$ and $n\geq \ee^{5b}$, $\overset{n-b}{\underset{k=1}{\prod}}\underset{u\in \Omega(w_k)}{\prod}p(u,x-L)\geq (1- \frac{\eta}{L^2})^\frac{1}{2} $. Plugging into (\ref{t12}) yields that
$$\Q((\xi_n)^c,w_n \text{ is a good vertex})\leq  \frac{\eta}{L^2} \Q(w_n\text{ is a a good vertex})\leq    \frac{\eta}{L^2} \Q(w_n\in \mathcal{Z}_n^{x,x+L,x-L}).$$
It follows from (\ref{t13}) that
$$\Q((\xi_n)^c, w_n\in \mathcal{Z}_n^{x,x+L,x-L})\leq  \eta(1+x)(\frac{1}{L^2}\Q(w_n\in Z^{x,x+L,x-L})+ n^{-\frac{3}{2}}).$$
Remember that the spine behaves as a centred random walk. Then apply (\ref{2.7}) to see that $\Q(w_n\in \mathcal{Z}_n^{x,x+L,x-L})\leq \alpha_{3}(1+L)^2n^{-\frac{3}{2}}$, it completes the proof of Lemma \ref{t17}.
\hfill $\Box$
\\

{\bf Acknowledgements:} I wish to thank my supervisor Yueyun Hu for introducing me to this subject and constantly finding the time for useful discussions and advice. I would also like to thank the referees for their helpful comments.

\bibliographystyle{plain}
\bibliography{bibli}

\end{document}